\documentclass[11pt]{amsart}
\usepackage{amsmath}
\usepackage{amsxtra}
\usepackage{amscd}
\usepackage{amsthm}
\usepackage{amsfonts}
\usepackage{amssymb}
\usepackage{eucal}

\theoremstyle{definition}

\theoremstyle{definition}
\newtheorem{thm}{Theorem}
\newtheorem{conj}{Conjecture}

\theoremstyle{remark}
\newtheorem{rem}{Remark}

\numberwithin{equation}{section}

\begin{document}

\newcommand{\thmref}[1]{Theorem~\ref{#1}}
\newcommand{\secref}[1]{Sect.~\ref{#1}}
\newcommand{\lemref}[1]{Lemma~\ref{#1}}
\newcommand{\propref}[1]{Proposition~\ref{#1}}
\newcommand{\corref}[1]{Corollary~\ref{#1}}
\newcommand{\remref}[1]{Remark~\ref{#1}}
\newcommand{\conjref}[1]{Conjecture~\ref{#1}}
\newcommand{\nc}{\newcommand}
\nc{\on}{\operatorname}
\nc{\ch}{\mbox{ch}}
\nc{\Z}{{\mathbb Z}}
\nc{\C}{{\mathbb C}}
\nc{\pone}{{\mathbb C}{\mathbb P}^1}
\nc{\pa}{\partial}
\nc{\F}{{\mathcal F}}
\nc{\arr}{\rightarrow}
\nc{\larr}{\longrightarrow}
\nc{\al}{\alpha}
\nc{\ri}{\rangle}
\nc{\lef}{\langle}
\nc{\W}{{\mathcal W}}
\nc{\la}{\lambda}
\nc{\ep}{\epsilon}
\nc{\su}{\widehat{{\mathfrak sl}}_2}
\nc{\sw}{{\mathfrak s}{\mathfrak l}}
\nc{\g}{{\mathfrak g}}
\nc{\h}{{\mathfrak h}}
\nc{\n}{{\mathfrak n}}
\nc{\N}{\widehat{\n}}
\nc{\G}{\widehat{\g}}
\nc{\De}{\Delta_+}
\nc{\gt}{\widetilde{\g}}
\nc{\Ga}{\Gamma}
\nc{\one}{{\mathbf 1}}
\nc{\z}{{\mathfrak Z}}
\nc{\zz}{{\mathcal Z}}
\nc{\Hh}{{\mathcal H}_\beta}
\nc{\qp}{q^{\frac{k}{2}}}
\nc{\qm}{q^{-\frac{k}{2}}}
\nc{\La}{\Lambda}
\nc{\wt}{\widetilde}
\nc{\qn}{\frac{[m]_q^2}{[2m]_q}}
\nc{\cri}{_{\on{cr}}}
\nc{\kk}{h^\vee}
\nc{\sun}{\widehat{\sw}_N}
\nc{\hh}{{\mathbf H}_{q,t}(\g)}
\nc{\HH}{{\mathcal H}_{q,t}}
\nc{\ca}{\wt{{\mathcal A}}_{h,k}(\sw_2)}
\nc{\si}{\sigma}
\nc{\gl}{\widehat{{\mathfrak g}{\mathfrak l}}_2}
\nc{\el}{\ell}
\nc{\s}{T}
\nc{\bi}{\bibitem}
\nc{\om}{\omega}
\nc{\WW}{\W_\beta}
\nc{\scr}{{\mathbf S}}
\nc{\ab}{{\mathbf a}}
\nc{\rr}{r}
\nc{\ol}{\overline}
\nc{\con}{qt^{-1} + q^{-1}t}
\nc{\den}{q^{\el-1} t^{-\el+1}+ q^{-\el+1} t^{\el-1}}
\nc{\ds}{\displaystyle}
\nc{\B}{B}
\nc{\A}{A^{(2)}_{2\el}}
\nc{\GG}{{\mathcal G}}
\nc{\UU}{{\mathcal U}}
\nc{\MM}{{\mathcal M}}
\nc{\CC}{{\mathcal C}}
\nc{\GL}{^L\G}
\nc{\dzz}{\frac{dz}{z}}
\nc{\Res}{\on{Res}}

\title[Deformations of ${\mathcal W}$--algebras]{Deformations of ${\mathcal
W}$--algebras associated to simple Lie algebras}

\author{Edward Frenkel}

\author{Nicolai Reshetikhin}

\address{Department of Mathematics, University of California, Berkeley, CA
94720, USA}

\date{July 1997; Revised September 1997}

\begin{abstract}
Deformed $\W$--algebra $\W_{q,t}(\g)$ associated to an arbitrary
simple Lie algebra $\g$ is defined together with its free field
realizations and the screening operators. Explicit formulas are given
for generators of $\W_{q,t}(\g)$ when $\g$ is of classical
type. These formulas exhibit a deep connection between $\W_{q,t}(\g)$
and the analytic Bethe Ansatz in integrable models associated to
quantum affine algebras $U_q(\G)$ and $U_t(\GL)$. The scaling limit of
$\W_{q,t}(\g)$ is closely related to affine Toda field theories.
\end{abstract}

\maketitle

\section{Introduction}

In this paper we define deformations of the $\W$--algebras associated
to arbitrary simple Lie algebras, their free field realizations and
the screening operators. The deformed $\W$--algebra $\W_{q,t}(\g)$
associated to a simple Lie algebra $\g$ is a family of associative
algebras depending on two parameters, $q$ and $t$. In the case when
$\g = A_\el$, the algebra $\W_{q,t}(\g)$ has been constructed in
\cite{SKAO,FF:w,AKOS} (see also \cite{LP1,FR:crit}). Various limits of
$\W_{q,t}(\g)$ can be identified with previously known algebras. In
particular, in the limit $q \arr 1$ with $t=q^\beta$ and $\beta$
fixed, we recover the ordinary $\W$--algebra corresponding to $\g$.

In the limits $t \arr 1$ and $q \arr 1$ the algebra $\W_{q,t}(\g)$
becomes commutative, and has a natural Poisson structure. We
conjecture that the Poisson algebra $\W_{q,1}(\g)$ is isomorphic to
the center of the quantized enveloping algebra $U_q(\G)$ at the
critical level (see \cite{FR:crit}), while the Poisson algebra
$\W_{1,t}(\g)$ is isomorphic to the Poisson algebra obtained by the
difference Drinfeld-Sokolov reduction of the loop group $G((z))$ (see
\cite{FRS,SS}). In the case $\g=A_\el$, the first conjecture follows
from results of \cite{FR:crit} and the second follows from results of
\cite{FRS,SS}; we have also checked the second conjecture for $\g=C_2$
(see Appendix B).

Another important limit is $q \arr \ep$, where $\ep=1$ for
simply-laced $\g$, and $\ep=\exp(\pi i/r^\vee)$ for nonsimply-laced
$\g$. Here $r^\vee$ is the order of the automorphism of a simply-laced
Lie algebra that gives rise to $\g$. In this limit, $\W_{q,t}(\g)$
contains a commutative subalgebra $\W'_{\ep,t}(\g)$, which has a
natural Poisson structure. We conjecture that the Poisson algebra
$\W'_{\ep,t}(\g)$ is isomorphic to the center of the quantized
enveloping algebra $U_t(^L\G)$ at the critical level. Here $^L\G$ is
the affine algebra that is Langlands dual to $\G$, that is the Cartan
matrix of $^L\G$ is the transpose of the Cartan matrix of $\G$.

We remark that for simply-laced $\g$, $\W_{q,t}(\g) \simeq
\W_{t,q}(\g)$, but for nonsimply-laced $\g$, $\W_{q,t}(\g)$ is not
isomorphic to $\W_{t,q}(^L\!\g)$, as one would expect by analogy with
Langlands duality of the ordinary $\W$--algebras \cite{FF:ds}. This
means that the duality becomes more complicated after deformation.

We give explicit formulas for generators of $\W_{q,t}(\g)$ when $\g$
is of classical type. These formulas exhibit a remarkable connection
with the analytic Bethe Ansatz in integrable models associated with
quantum affine algebras (see \cite{Ba,Re1,Re2,BR,KS1}). Recall that in
\cite{FR:crit} we conjectured (and proved in the case $\g=A_\el$) that
the formulas for the free field realization of the center of
$U_q(\G)$, i.e., $\W_{q,1}(\g)$, coincide with Bethe Ansatz formulas
for the eigenvalues of the transfer-matrices in the $U_q(\G)$
integrable model. In this paper we will see further evidence of
that. We will also see that the formulas for the free field
realization of $\W'_{\ep,t}(\g)$ coincide with Bethe Ansatz formulas
for transfer-matrices in the $U_t(^L\G)$ model. Furthermore, the
formulas for the free field realization of $\W_{1,t}(\g)$ exhibit
similarity with formulas for the eigenvalues of transfer-matrices in
the integrable model associated to $U_t(\G^\vee)$, where
$\G^\vee = \; ^L\!(\widehat{^L\!\g})$.

Thus, the free field realization of $\W_{q,t}(\g)$ connects the
analytic Bethe Ansatz formulas for $U_q(\G)$, $U_t(^L\G)$ and
$U_t(\G^\vee)$.

The Bethe Ansatz formulas can be written for any finite-dimensional
representation of quantum affine algebra, and can be thought of as
$q$--analogues of the character formulas. Numerous examples are known
in the literature, see \cite{Re2,KS1,KS2}. However, there seems to be
no systematic algorithm for writing these formulas in general. We hope
that one can use the property of commutativity with the screening
operators as the defining property for these ``$q$--characters''. This
may give us some insights into the category of finite-dimensional
representations of quantum affine algebras. Furthermore, it would then
appear that $\W_{q,t}(\g)$ is a simultaneous deformation of
representation rings of $U_q(\G)$ and $U_t(^L\G)$. We will discuss
these issues in more detail in \cite{FR:new}.

The paper is organized as follows: in Sect.~2 we define a
two-parameter family of deformations of the Cartan matrix associated
to each simple Lie algebra. In Sect.~3 we define the Heisenberg
algebra ${\mathcal H}_{q,t}(\g)$ and the screening operators. We then
define $\W_{q,t}(\g)$ as the centralizsr of the screening operators in
${\mathcal H}_{q,t}(\g)$. We conjecture the form of the generators of
$\W_{q,t}(\g)$ and derive from this conjecture the exchange relations
between the generators of $\W_{q,t}(\g)$.

We also compute the relations between the screening currents. It was
shown in \cite{FF:w} that for $\g=A_\el$ they satisfy elliptic
relations, which can be considered as elliptic analogues of Drinfeld's
relations in $U_q(\widehat{{\mathfrak n}})$ \cite{Dr}. This fact has
been used in \cite{Mi1,Mi2} to construct two-sided resolutions of
irreducible representations of deformed $\W$--algebras. Here we obtain
the relations for arbitrary $\g$.

In Sect.~4 we study various limits of $\W_{q,t}(\g)$ and identify them
with some known algebras. In Sect.~5 we give explicit formulas for the
generating field $T_1(z)$ of $\W_{q,t}(\g)$, in the case when $\g$ is
a simple Lie algebra of classical type. In Sect.~6 we discuss the
connection between $\W_{q,t}(\g)$ and the analytic Bethe Ansatz. In
Sect.~7 we define the deformed $\W$--algebra associated to the
self-dual twisted affine algebra $A^{(2)}_{2\el}$ generalizing the
recent work \cite{BL}.

In Sect.~8, we consider the scaling limit of the exchange relations
between the generators of $\W_{q,t}(\g)$ and between the screening
currents. We show that the scaling limit of $\W_{q,t}(\g)$ can be
identified with the Faddeev-Zamolodchikov algebra of an affine Toda
field theory. In the case $\g = A_\el$ this has been suggested by
S.~Lukyanov \cite{Lu2,Lu3} (see also \cite{Lu1}). The scaling limit of
our free field realization can therefore be used to obtain explicit
formulas for form-factors in general affine Toda field theories along
the lines of \cite{Lu2,Lu3,BL}. On the other hand, relations between
the screening currents give rise in the scaling limit to Drinfeld's
relations, for both twisted and non-twisted affine algebras.

In Sect.~9 we define the analogues of vertex operators (primary
fields) for $\W_{q,t}(\g)$, which correspond to the fundamental
representations of $U_q(\G)$. Using these vertex operators and the
screening operators, one can give complete bosonization of the vertex
operators in the general solvable SOS models, along the lines of
\cite{LP2,AJMP}. We note that the last two sections are somewhat
more physics oriented than the rest of the paper.

In Appendix A we give some explicit formulas for the Poisson algebras
${\mathcal H}_{q,1}(\g)$ and $\W_{q,1}(\g)$ for Lie algebras of
classical types. In Appendix B we discuss the difference
Drinfeld-Sokolov reduction in the case $\g=C_2$. Appendix C contains
explicit formulas for the matrices $M(q,t)$ corresponding to the
classical Lie algebras. In Appendix D we recall the definition of
deformed chiral algebra from \cite{FR:dca}.

\medskip
\noindent{\bf Acknowledgements.} The research of the first author was
supported by grants from the Packard Foundation and the NSF. The
research of the second author was supported by an NSF grant DMS
9401163.

\section{Two-parameter deformations of Cartan matrices}

\subsection{General formula}

Let $\g$ be a simple Lie algebra of rank $\el$. Let $(\cdot,\cdot)$ be the
invariant inner product on $\g$, normalized as in \cite{Kac}, so that the
square of the maximal root equals $2$. Let $\{ \al_1,\ldots,\al_\el \}$ and
$\{ \om_1,\ldots,\om_\el \}$ be the sets of simple roots and of fundamental
weights of $\g$, respectively. We have:
$$
(\al_i,\om_j) = \frac{(\al_i,\al_i)}{2} \delta_{i,j}.
$$
Let $r^\vee$ be the maximal number of edges connecting two vertices of the
Dynkin diagram of $\g$. Thus, $r^\vee=1$ for simply-laced $\g$, $r^\vee=2$
for $B_\el, C_\el, F_4, G_2$, and $r^\vee=3$ for $D_4$. Set
$$
D = \on{diag}(\rr_1,\ldots,\rr_\el),
$$
where
\begin{equation}    \label{di}
\rr_i = r^\vee \frac{(\al_i,\al_i)}{2}.
\end{equation}
All $\rr_i$'s are integers; for simply-laced $\g$, $D$ is the identity
matrix.

Now let $C = (C_{ij})_{1\leq i,j\leq \el}$ be the {\em Cartan matrix} of
$\g$. We have:
$$
C_{ij} = \frac{2(\al_i,\al_j)}{(\al_i,\al_i)}.
$$
Denote by $I = (I_{ij})_{1\leq i,j\leq \el}$ the {\em incidence matrix},
$$
I_{ij} = 2 \delta_{i,j} - C_{ij}.
$$
Let $\B = (\B_{ij})_{1\leq i,j\leq \el}$ be the following matrix:
$$
B = D C,
$$
i.e.,
$$
\B_{ij} = r^\vee (\al_i,\al_j).
$$

Now let $q, t$ be indeterminates. We will use the standard notation
$$
[n]_q = \frac{q^n - q^{-n}}{q - q^{-1}}.
$$
We define $\el \times \el$ matrices $C(q,t)$, $D(q,t)$, and $\B(q,t)$ by the
formulas
\begin{align}    \label{qtc}
C_{ij}(q,t) &= (q^{\rr_i} t^{-1} + q^{-\rr_i} t) \delta_{i,j} - [I_{ij}]_q,
\\ D(q,t) &= \on{diag}([\rr_1]_q,\ldots,[\rr_\el]_q), \label{qtd} \\
\B(q,t) &= D(q,t) C(q,t). \notag
\end{align}
Thus,
\begin{equation}    \label{qts}
\B_{ij}(q,t) = [\rr_i]_q \left( (q^{\rr_i} t^{-1} + q^{-\rr_i} t)
\delta_{i,j} - [I_{ij}]_q \right).
\end{equation}
It is easy to see that the matrix $\B(q,t)$ is symmetric. For simply-laced
$\g$,
$$
C_{ij}(q,t) = \B_{ij}(q,t) = (q t^{-1} + q^{-1} t) \delta_{i,j} - I_{ij}.
$$

Clearly, the limits of $C(q,t)$, $D(q,t)$, and $\B(q,t)$ as $q \arr 1$ and
$t \arr 1$ coincide with $C$, $D$, and $\B$, respectively. Note also that
$$
\B_{ij}(q,1) = [\B_{ij}]_q,
$$
and
$$
\B_{ij}(1,t) = \rr_i ((t + t^{-1}) \delta_{ij} - I_{ij}).
$$

It is interesting that for all $\g$, the Coxeter number $h$ enters the
determinant of $C(q,t)$ as the power of $t$, and $r^\vee h^\vee$,
where $h^\vee$ is the dual Coxeter number of $\g$, enters the
determinant as the power of $q$. For instance, for $\g=A_\el$:
$$
\on{det} C(q,t) = \frac{q^{\el+1} t^{-\el-1} - q^{-\el-1} t^{\el+1}}{q t^{-1}
- q^{-1} t}
$$
(in this case $h = r^\vee h^\vee=\el+1$), for $\g=B_\el$:
$$
\on{det} C(q,t) = q^{2\el-1} t^{-\el} + q^{-2\el+1} t^{\el}
$$
(in this case $h=2\el, r^\vee h^\vee = 2(2\el-1)$), for $\g=C_\el$:
$$
\on{det} C(q,t) = q^{\el+1} t^{-\el} + q^{-\el-1} t^{\el}
$$
(in this case $h=2\el, r^\vee h^\vee = 2(\el+1)$), etc.

\section{Deformed $\W$--algebras}

\subsection{Heisenberg algebra $\HH(\g)$}

Let $\HH(\g)$ be the Heisenberg algebra with generators $a_i[n],
i=1,\ldots,\el; n \in \Z$, and relations
\begin{equation}    \label{a}
[a_i[n],a_j[m]] = \frac{1}{n} (q^n - q^{-n}) (t^n - t^{-n}) \B_{ij}(q^n,t^n)
\delta_{n,-m}.
\end{equation}
In this and other formulas of this type we understand that the $0$th
generator commutes with all other generators: $[a_i[0],a_j[m]] = 0, \forall
m \in \Z$.

The generators $a_i[n]$ are ``root'' type generators of $\HH(\g)$. There is
a unique set of ``fundamental weight'' type generators, $y_i[n],
i=1,\ldots,\el; n \in \Z$, that satisfy:
\begin{equation}    \label{ay}
[a_i[n],y_j[m]] = \frac{1}{n} (q^{\rr_i n} - q^{-\rr_i n})(t^n - t^{-n})
\delta_{i,j} \delta_{n,-m}.
\end{equation}
They have the following commutation relations:
\begin{equation}    \label{y}
[y_i[n],y_j[m]] = \frac{1}{n} (q^{n} - q^{-n}) (t^n - t^{-n})
M_{ij}(q^n,t^n) \delta_{n,-m},
\end{equation}
where $(M_{ij}(q,t))_{1\leq i,j\leq \el}$ is the following matrix
\begin{align}    \label{tildeM}
M(q,t) &= D(q,t) C(q,t)^{-1} \\
&= D(q,t) B(q,t)^{-1} D(q,t).  \notag
\end{align}

Note that
\begin{equation}    \label{expr}
a_j[n] = \sum_{j=1}^\el C_{ij}(q^n,t^n) y_j[n].
\end{equation}

Introduce the generating series:
$$
A_i(z) = t^{2(\rho^\vee,\al_i)} q^{-2r^\vee (\rho,\al_i) + 2a_i[0]}
:\exp \left( \sum_{m\neq 0} a_i[m] z^{-m} \right):,
$$
$$
Y_i(z) = t^{2(\rho^\vee,\om_i)} q^{-2r^\vee (\rho,\om_i) + 2 y_i[0]}
:\exp \left( \sum_{m\neq 0} y_i[m] z^{-m} \right):.
$$
Note that $(\rho^\vee,\al_i) = 1, r^\vee(\rho,\al_i) = \rr_i$.

For each weight $\mu$ of the Cartan subalgebra of $\g$, let $\pi_\mu$ be
the Fock representation of $\HH(\g)$ generated by a vector $v_\mu$, such
that $a_i[n] v_\mu = 0, n>0$, and $a_i[0] v_\mu = (\mu,\al_i) v_\mu$.

\subsection{Screening operators}

Set $t = q^\beta$. Introduce the shift operators $e^{Q_i}, i=1,\ldots,\el$,
acting from $\pi_\mu$ to $\pi_{\mu+\beta\al_i}$, which satisfy commutation
relations
\begin{equation}    \label{qj}
[a_i[n],e^{Q_j}] = \B_{ij} \beta \delta_{n,0} e^{Q_j}.
\end{equation}

Let
\begin{equation}    \label{sm+}
s^+_i[m] = \frac{a_i[m]}{q^{m \rr_i}-q^{-m \rr_i}}, \quad m\neq 0, \quad \quad
s^+_i[0] = a_i[0]/\rr_i,
\end{equation}
\begin{equation}    \label{sm-}
s^-_i[m] = \frac{a_i[m]}{t^m-t^{-m}}, \quad m\neq 0,
\quad \quad s^-_i[0] = a_i[0]/\beta.
\end{equation}

Now define the {\em screening currents} by the formulas
\begin{align}    \label{s+}
S_i^+(z) & = e^{-Q_i/\rr_i} z^{-s^+_i[0]} :\exp \left( \sum_{m\neq 0}
s^+_i[m] z^{-m} \right):,\\    \label{s-}
S_i^-(z) & = e^{Q_i/\beta} z^{s^-_i[0]} :\exp \left( - \sum_{m\neq 0}
s^-_i[m] z^{-m} \right):.
\end{align}

They satisfy the difference equations:
\begin{equation}    \label{scr1}
S^+_i(zq^{-\rr_i}) = t^{-2} q^{2\rr_i} :A_i(z) S^+_i(zq^{\rr_i}):,
\end{equation}
and
\begin{equation}    \label{scr2}
S^-_i(zt) = t^{-2} q^{2\rr_i} :A_i(z) S^-_i(zt^{-1}):.
\end{equation}

\medskip
\noindent{\bf Remark on notation.} To avoid confusion, let us
emphasize that in the case $\g=A_\el$ our notation here differs
slightly from that of \cite{FF:w}. Namely, our $q$ and $t$ here
correspond to $q^{1/2}$ and $(q/p)^{1/2}$, respectively, of
\cite{FF:w} (though our $\beta$ coincides with $\beta$ of
\cite{FF:w}). We made this change of notation to avoid the appearance
of half-integers in the formulas. Also, this notation agrees with that
of \cite{FR:crit}.

The connection between our notation and that of \cite{AJMP} is as follows:
$x=q/t, r=1/(1-\beta)$; $q=x^r, t=x^{r-1}$.

\subsection{Definition of $\W_{q,t}(\g)$}

Let ${\mathbf H}_{q,t}(\g)$ be the vector space spanned by formal power
series of the form
\begin{equation}    \label{mono}
:\pa_z^{n_1} Y_{i_1}(zq^{j_1}t^{k_1})^{\ep_1} \ldots \pa_z^{n_1}
Y_{i_l}(zq^{j_l}t^{k_l})^{\ep_l}:,
\end{equation}
where $\ep_i = \pm 1$. The pair $(\hh,\pi_0)$ is a {\em deformed
chiral algebra} (DCA) in the sense of \cite{FR:dca}. The definition of
DCA is recalled in Appendix D.

Denote
\begin{align*}
S^+_i &= \int S^+_i(z) dz: \pi_0 \arr \pi_{-\beta \al_i/\rr_i}, \\
S^-_i &= \int S^-_i(z) dz: \pi_0 \arr \pi_{\al_i}.
\end{align*}
Here the integral simply means the $(-1)$st Fourier coefficient of the
series $S^\pm_i(z)$, which is a well-defined linear operator. We call
these operators the {\em screening operators}.

We define the DCA ${\mathbf W}_{q,t}(\g)$ as the maximal subalgebra of
$(\hh,\pi_0)$, which commutes with the operators $S^\pm_i, i=1,\ldots,\el$,
i.e., the subspace of $\hh$, which consists of all fields that commute with
these operators. We define the deformed $\W$--algebra $\W_{q,t}(\g)$ as the
associative algebra, topologically generated by the Fourier coefficients of
fields from ${\mathbf W}_{q,t}(\g)$. All elements of the algebra
$\W_{q,t}(\g)$ act on the Fock representations $\pi_\la$ and commute with
the screening operators.

We call the field \eqref{mono} {\em elementary}, if it does not contain
derivatives. We assign to such a term the element of the weight lattice of
$\g$,
$$
\sum_{a=1}^\el \ep_a \omega_{i_a}.
$$

\begin{conj}    \label{conj1}
{\em Let $\g$ be a simple Lie algebra and $\G$ be the corresponding
affine Kac-Moody algebra. For each $i=1,\ldots,\el$, there exists a
field $T_i(z)$ in ${\mathbf W}_{q,t}(\g)$, such that $T_i(z) = Y_i(z)
+$ the sum of elementary terms of the form
$$
c(q,t) :Y_i(z) A_{i_1}(zq^{a_1} t^{b_1})^{-1} \ldots A_{i_k}(zq^{a_k}
t^{b_k})^{-1}:
$$
(where $c(q,1)$ is a positive integer independent of
$q$). Furthermore, the set of weights of these terms counted with
multiplicity $c(q,1)$ is the set of weights of the finite-dimensional
irreducible representation $V_{\omega_i}$ of $U_q(\G)$ with highest
weight $\omega_i$.}
\end{conj}

Such fields $T_i(z)$ have been constructed in \cite{SKAO,FF:w,AKOS} in
the $A_\el$ case. In \secref{class-types} we will explicitly construct
the field $T_1(z)$ for all simple Lie algebras of classical
types. Some motivations for \conjref{conj1} are discussed in
Sect.~6.1.

\subsection{Exchange relations}

In this subsection we compute the exchange relations between the fields
$T_i(z)$, provided that our \conjref{conj1} is true. We assume
throughout the rest of this section that $|q|<1$ and $|t|<1$.

According to \conjref{conj1}, each $T_i(z)$ is the sum of $Y_i(z)$ and
other terms which are normally ordered products of $Y_i(z)$ and
$A_i(zq^a t^b)^{-1}$. Let us recall from \cite{FR:dca} what we mean by
an exchange relation. It follows from the commutation relations
\eqref{a}--\eqref{y} that the composition $T_i(z) T_j(w)$ converges
for $|z| \gg |w|$ and can be analytically continued to a meromorphic
function on $\C^\times \times \C^\times$ with poles on shifted
diagonals $z=w\gamma$, where $\gamma \in q^{\Z} t^{\Z}$. Let us denote
this analytic continuation by $R(T_i(z) T_j(w))$. By exchange relation
we understand a relation of the type
$$
R(T_i(z) T_j(w)) = S_{T_i,T_j} \left( \frac{w}{z} \right) R(T_j(w) T_i(z)),
$$
where $S(z)$ is a meromorphic function.

According to formulas \eqref{a} and \eqref{ay} the analytic continuations of
$Y_i(z)$ and $A_i(w)^{-1}$ satisfy:
\begin{align}    \label{local}
R(Y_i(z) A_j(w)^{-1}) &= R(A_j(w)^{-1} Y_i(z)), \\ \notag R(A_i(z)^{-1}
A_j(w)^{-1}) &= R(A_j(w)^{-1} A_i(z)^{-1}),
\end{align}
so that fields $A_j(z)^{-1}$ are mutually local and also local with the
fields $Y_i(w)$. Formulas \eqref{local} imply that if $G^\al_i(z)$ is one
of the terms entering $T_i(z)$ and $G^\beta_j(w)$ is one of the terms
entering $T_j(w)$, then
$$
S_{G^\al_i,G^\beta_j} \left( \frac{w}{z} \right) = S_{Y_i,Y_j} \left(
\frac{w}{z} \right),
$$
which means that
$$
S_{T_i,T_j} \left( \frac{w}{z} \right) = S_{Y_i,Y_j} \left(
\frac{w}{z} \right).
$$

It is straightforward to find from \eqref{y} that
\begin{equation}    \label{theta}
S_{Y_i,Y_j} \left( \frac{w}{z} \right) = \exp \sum_{n>0} \frac{1}{n}
(q^n-q^{-n})(t^n-t^{-n})M_{ij}(q^n,t^n) \left( \left( \frac{w}{z}
\right)^n - \left( \frac{z}{w} \right)^n \right),
\end{equation}
where the matrix $(M_{ij})_{1\leq i,j\leq \el}$ is given by formula
\eqref{tildeM}. Formula \eqref{theta} is an elliptic function with the
multiplicative period $q^{2r^\vee h^\vee} t^{2h}$, where $h$ and $h^\vee$
are the Coxeter and the dual Coxeter numbers, respectively. It can be
rewritten as the ratio of products of theta-functions, as in \cite{FF:w} in
the case of $A_\el$.

One can also obtain quadratic relations on the Fourier coefficients of
various fields from the DCA ${\mathbf W}_{q,t}(\g)$ similar to the
relations found in \cite{SKAO,FF:w,AKOS}. A general method for
writing such relations is described in \cite{FR:dca}.

\subsection{Relations between the screening currents}

Let $p_i = q^{\rr_i} t^{-1}$, $\rr_{ij} = \on{min}(\rr_i,\rr_j)$, and
$p_{ij}=q^{\rr_{ij}} t^{-1}$. Denote
$$
^L\!\B_{ij} = \frac{B_{ij}}{\rr_i\rr_j}.
$$
Introduce the notation
$$
\theta(z;a) = \prod_{n=1}^\infty (1-zq^{n-1})(1-z^{-1}q^n)(1-q^n).
$$

Direct computation gives us the following relations on the analytic
continuations of the compositions of the screening currents $S^+_i(z)$:
\begin{align*}
S^+_i(z) S^+_i(w) &= p_i^{-2} \left( \frac{w}{z} \right)^{-^L\!\B_{ii}\beta+2}
\frac{\theta \left( \ds \frac{w}{z} p_i^2 ; q^{2\rr_{ii}}
\right)}{\theta \left( \ds \frac{w}{z} p_i^{-2} ; q^{2\rr_{ii}}
\right)} S^+_i(w) S^+_i(z), \\
S^+_i(z) S^+_j(w) &= - p_{ij} \left( \frac{w}{z} \right)^{-^L\!\B_{ij}\beta-1}
\frac{\theta \left( \ds \frac{w}{z} p_{ij}^{-1} ; q^{2\rr_{ij}}
\right)}{\theta \left( \ds \frac{w}{z} p_{ij} ;
q^{2\rr_{ij}} \right)} S^+_j(w) S^+_i(z), \quad i\neq j,
\B_{ij} \neq 0, \\
S^+_i(z) S^+_j(w) &= S^+_j(w) S^+_i(z), \quad \quad \B_{ij}=0.
\end{align*}

Let $\wt{p}_{ij} = q^{\B_{ij}} t^{-1}$. We have the following relations for
$S^-_i(z)$:
\begin{align*}
S^-_i(z) S^-_i(w) &= p_i^2 \left( \frac{w}{z}
\right)^{-\B_{ii}/\beta+2} \frac{\theta \left( \ds \frac{w}{z} p_i^{-2} ;
t^2 \right)}{\theta \left( \ds \frac{w}{z} p_i^2 ; t^2 \right)}
S^-_i(w) S^-_i(z), \\ S^-_i(z) S^-_j(w) &= - \wt{p}_{ij}^{-1} \left(
\frac{w}{z} \right)^{-\B_{ij}/\beta-1} \frac{\theta \left( \ds \frac{w}{z}
\wt{p}_{ij} ; t^2 \right)}{\theta \left( \ds \frac{w}{z}
\wt{p}_{ij}^{-1} ; t^2 \right)} S^-_j(w) S^-_i(z), \quad \quad i\neq j.
\end{align*}

For $\g=A_\el$ these relations were obtained in \cite{FF:w}.

One can use the screening currents to construct resolutions of
irreducible representations of $\W(\g)$, as in \cite{LP2,Mi1,Mi2} for
$\g=A_\el$. For that one needs to multiply each screening current by a
ratio of theta-functions, to make them single-valued, as in
\cite{LP2,Mi1,Mi2}.

\section{Identification of various limits of $\W_{q,t}(\g)$}

\subsection{The conformal limit $q \arr 1, t = q^\beta$}

Let $\Hh(\g)$ be the Heisenberg algebra with generators $\ab_i[n],
i=1,\ldots,\el; n \in \Z$, and relations:
\begin{equation}    \label{oh}
[\ab_i[n],\ab_j[m]] = n \B_{ij} \beta \delta_{n,-m}.
\end{equation}
By abuse of notation, we denote by $\pi_\mu$ the Fock representation of
$\Hh(\g)$ generated by a vector $v_\mu$, such that $\ab_i[n] v_\mu = 0,
n>0$, and $\ab_i[0] v_\mu = (\mu,\al_i) v_\mu$. We also use the same
notation as before, $e^{Q_i}$, for the shift operators acting from
$\pi_\mu$ to $\pi_{\mu+\beta\al_i}$ and satisfying
$$
[\ab_i[n],e^{Q_j}] = \B_{ij} \beta \delta_{n,0} e^{Q_j}.
$$

Now set
\begin{align}    \label{s+c}
\scr_i^+(z) & = e^{-Q_i/\rr_i} z^{-\ab_i[0]/\rr_i} :\exp \left( \sum_{m\neq 0}
\frac{1}{m\rr_i} \ab_i[m] z^{-m} \right):,\\ \label{s-c} \scr_i^-(z) & =
e^{Q_i/\beta} z^{\ab_i[0]/\beta} :\exp \left( - \sum_{m\neq 0}
\frac{1}{m\beta} \ab_i[m] z^{-m} \right):.
\end{align}
Denote
$$
\scr_i^\pm = \int \scr_i^\pm(z) dz.
$$
These are the ordinary screening operators.

Let us recall the definition of the ordinary $\W$--algebra associated to
$\g$ from \cite{FF:ds,FF:laws}. Let $\pi_0$ be the vertex operator algebra
(VOA) associated to the Heisenberg algebra $\Hh(\g)$ (see
\cite{FF:laws}). For generic $\beta$ the VOA ${\mathbf W}_\beta(\g)$ is
defined as the vertex operator subalgebra of the VOA $\pi_0$, which is the
intersection of kernels of the screening operators $\scr_i^-,
i=1,\ldots\el$:
$$
{\mathbf W}_\beta(\g) = \bigcap_{i=1}^\el \on{Ker} \scr_i^-.
$$
Thus, ${\mathbf W}_\beta(\g)$ consists of the fields that commute with
$\scr_i^-$. The $\W$--algebra $\WW(\g)$ is defined as the associative (or
Lie) algebra generated by the Fourier coefficients of the fields from
${\mathbf W}_\beta(\g)$.

\begin{rem}
In the notation of \cite{FF:laws}, our $\WW(\g)$ is $\W_\gamma(\g)$, where
$\gamma = (r^\vee/\beta)^{1/2}$.\qed
\end{rem}
\medskip

It was proved in \cite{FF:ds} (see also \cite{FF:laws}) that for generic
$\beta$, ${\mathbf W}_\beta(\g)$ automatically lies in the kernel of the
other set of screening operators, $\scr_i^+, i=1,\ldots,\el$. This implies
the following duality.

\begin{thm}[\cite{FF:ds}]    \label{duality}
{\em For generic $\beta$,
$$
\WW(\g) \simeq \W_{r^\vee/\beta}(^L\!\g),
$$
where $^L\!\g$ is the Langlands dual Lie algebra to $\g$.}
\end{thm}

Now let us consider the limit of $\W_{q,t}(\g)$ as $q \arr 1$ with $t =
q^\beta$. Let us pass to a new set of generators
$$
\wt{a}_i[n] = \frac{a_i[n]}{q-q^{-1}}.
$$
It is clear from formula \eqref{a} that in the limit the operators
$\wt{a}_i[n]$ satisfy the relations given by formula \eqref{oh}. Hence
we can identify $\wt{a}_i[n]$ with $\ab_i[n]$ in this
limit. Furthermore, it is easy to see that the limit of the field
$S_i^\pm(z)$ equals $\scr_i^\pm(z)$. Comparing the definitions of
$\W_{q,t}(\g)$ and $\W_\beta(\g)$, we obtain

\begin{thm}
{\em In the limit $q \arr 1$ with $t = q^\beta$, $\W_{q,t}(\g)$ becomes
$\W_\beta(\g)$.}
\end{thm}

It is clear from the definition that $\W_{q,t}(\g) \simeq
\W_{q^{-1},t^{-1}}(\g)$ for all $\g$, and for simply-laced $\g$,
$$
\W_{q,t}(\g) \simeq \W_{t,q}(\g).
$$
The latter is the analogue of \thmref{duality} (note that $^L\!\g = \g$ in
the simply-laced case). However, for nonsimply-laced $\g$, $\W_{q,t}(\g)$
is {\em not} isomorphic to $\W_{t,q}(^L\!\g)$.

\subsection{The first classical limit $t \arr 1$}

Let us consider the limit $\W_{q,1}(\g)$ of $\W_{q,t}(\g)$ as $t \arr 1$
($q$ is fixed). Both ${\mathcal H}_{q,1}(\g)$ and $\W_{q,1}(\g)$ are
commutative algebras, but they inherit a Poisson structure. The
Poisson-Heisenberg algebra ${\mathcal H}_{q,1}(\g)$ has the relations
$$
\{ a_i[n],a_j[n] \} = (q^{\B_{ij}n} - q^{-\B_{ij}n}) \delta_{n,-m}.
$$
The Poisson algebra ${\mathcal H}_{q,1}$ was introduced in
\cite{FR:crit}. $\W_{q,1}(\g)$ is its Poisson subalgebra. It is easy to see
that the limit $q \arr 1$ of $\W_{q,1}(\g)$ is isomorphic to the limit of
$\WW(\g)$ as $\beta \arr 0$. The latter is known to be isomorphic to the
center of $U(\widehat{\g})$ at the critical level and to $\W(^L\!\g)$
\cite{FF:ds}.

\begin{conj}    \label{center}
{\em The limit $t \arr 1$ of $\W_{q,t}(\g)$ with fixed $q$ is
isomorphic, as a Poisson algebra, to the center of $U_q(\widehat{\g})$
at the critical level.}
\end{conj}

\conjref{center} was proved in \cite{FR:crit} for $\g = A_\el$. It has
been verified in \cite{FR:BC} for $\g = B_\el, C_\el$ and in \cite{Ko}
for $\g = D_\el, E_6$ and $G_2$. Some of these results are given in
Appendix A.

The following observation also confirms the conjecture. Consider the
quantized enveloping algebra $U_q(\G)$ and its loop-like generators
(Drinfeld's new generators), $\kappa_{i,n}$ in the notation of
\cite{Dr}. In the limit when the level tends to $0$, they generate a
Poisson algebra ${\mathcal B}_q(\g)$. In \cite{FR:crit} we explained
that the free field realization of $U_q(\G)$ induces an embedding of
the center of $U_q(\G)$ into a Heisenberg-Poisson algebra ${\mathcal
A}_q(\g)$, which is part of the free field realization of
$U_q(\G)$. In the case $\g=A_\el$, when the free field realization is
available \cite{AOS}, we know that this Heisenberg-Poisson algebra
${\mathcal A}_q(\g)$ is isomorphic to ${\mathcal B}_q(\g)$. We expect
that the same is true for other $\g$. But it is easy to see that
${\mathcal B}_q(\g)$ is isomorphic to ${\mathcal H}_{q,1}(\g)$.

\subsection{The second classical limit $q \arr 1$}

Now consider the limit $\W_{1,t}(\g)$ of $\W_{q,t}(\g)$ as $q \arr 1$
($t$ is fixed). This is again a Poisson algebra, which is a Poisson
subalgebra of ${\mathcal H}_{1,t}(\g)$. The latter has the following
relations:
$$
\{ a_i[n],a_j[n] \} = \rr_i ((t + t^{-1}) \delta_{ij} - I_{ij})
\delta_{n,-m}.
$$
The $t \arr 1$ limit of $\W_{1,t}(\g)$ coincides with the $\beta \arr
\infty$ limit of $\WW(\g)$, which is isomorphic to the classical
$\W$--algebra $\W(\g)$ obtained by the Drinfeld-Sokolov reduction of
the dual space to $\widehat{\g}$ (see \cite{FF:ds}). On the other
hand, in \cite{FRS,SS} a $p$--difference analogue of the
Drinfeld-Sokolov reduction was defined as a Poisson reduction of the
Poisson-Lie group $\GG=G((z))$. The result of this reduction is a
$p$--deformation of $\W(\g)$, which we denote by $\W^p(\g)$. There is
also an embedding of $\W^p(\g)$ into a Heisenberg-Poisson algebra;
this is a difference analogue of the Miura transformation from
\cite{DS}. These constructions are recalled in Appendix B.

\begin{conj}    \label{tds}
{\em The limit $q \arr 1$ of $\W_{q,t}(\g)$ with fixed $t$ is
isomorphic to $\W^p(\g)$, where $p=t^{r^\vee}$. Furthermore, the free
field realization of $\W_{1,t}(\g)$ coincides with the difference
Miura transformation of $\W^p(\g)$.}
\end{conj}

For simply-laced $\g$, the two Poisson algebras obtained as limits $q
\arr 1$ and $t \arr 1$ of $W_{q,t}(\g)$ are isomorphic: $\W_{q,1}(\g)
\simeq \W_{1,q}(\g)$. On the other hand, for $\g = A_\el$, it follows
from \cite{FRS,SS} that $\W_{t,1}(A_\el)$ is the Poisson algebra
obtained by the $t$--difference Drinfeld-Sokolov reduction of
$\widehat{SL}_{\el+1}$. Hence \conjref{tds} holds for $\g=A_\el$.

In Appendix B we check \conjref{tds} for $\g=C_2$.

\begin{rem}
We know that the center of $U(\widehat{\g})$ at the critical level is
isomorphic to $\W(^L\!\g)$. However, for nonsimply-laced $\g$ and $q
\neq 1$, the center of $U_q(\widehat{\g})$ is not isomorphic to
$\W_{1,q}(^L\!\g)$.\qed
\end{rem}

\subsection{The limit $q \arr \ep$}

Let $\ep=1$ for simply-laced $\g$, and $\ep=\exp(\pi i/r^\vee)$ for
nonsimply-laced $\g$. Since we have already discussed the limit $q
\arr 1$, we can now focus on nonsimply-laced $\g$. Let $R_s$
(respectively, $R_l$) be the subset of $\{ 1,2,\ldots,\el \}$,
consisting of the labels of short (respectively, long) simple
roots. Note that $r_i = 2$ for $i \in R_s$, and $r_i = 2r^\vee$ for $i
\in R_l$. The inspection of formulas \eqref{a}, \eqref{ay}, and
\eqref{y} shows that although the algebra ${\mathcal H}_{q,t}(\g)$
is not commutative when $q=\ep$, it contains a large center ${\mathcal
H}'_t(\g)$, generated by $a_i[n], i \in R_l, n \in \Z$, and
$a_i[r^\vee n], i \in R_s, n \in \Z$. The center ${\mathcal H}'_t(\g)$
has a natural Poisson structure, with the Poisson brackets between the
generators given by
\begin{align*}
\{ a_i[n],a_i[m] \} & = (-1)^{n+1} (t^{2n} - t^{-2n}) \delta_{n,-m}, & i \in
R_l, \\
\{ a_i[r^\vee n],a_i[r^\vee m] \} & = (-1)^{n+1} \frac{1}{r^\vee}
(t^{2 r^\vee n} - t^{-2 r^\vee n}) \delta_{n,-m}, & i \in
R_s, \\
\{ a_i[n],a_j[m] \} & = (-1)^n (t^{n} - t^{-n})
\delta_{n,-m}, & i,j \in R_l; I_{ij} \neq 0, \\ 
\{ a_i[r^\vee n],a_j[m] \} & = (-1)^n (t^{r^\vee n} - t^{- r^\vee
n}) \delta_{r^\vee n,-m}, & i \in R_s, j \in R_l; I_{ij}
\neq 0, \\
\{ a_i[r^\vee n],a_j[r^\vee m] \} & = (-1)^n \frac{1}{r^\vee} (t^{r^\vee
n} - t^{- r^\vee n}) \delta_{n,-m}, & i,j \in R_s; I_{ij}
\neq 0.
\end{align*}

Now the Dynkin diagram of $\g$ is obtained by folding from the Dynkin
diagram of a simply-laced Lie algebra $\wt{\g}$ under the action of an
automorphism of order $r^\vee$. Let $\Gamma$ be the set of vertices of
the Dynkin diagram of $\wt{\g}$, and $\sigma: \Gamma \arr \Gamma$ be
the corresponding automorphism of order $r^\vee$. Denote by
$(\wt{B}_{ij})$ the Cartan matrix of $\wt{\g}$. Consider the Poisson
algebra $\wt{{\mathcal H}}_t(\g)$ with generators $b_i[n], i \in
\Gamma, n \in \Z$, and relations:
\begin{align}    \label{dr1}
\{ b_i[n],b_j[m] \} &= \sum_{k=0}^{r^\vee} (t^{n\wt{B}_{i,\sigma(j)}} -
t^{-n\wt{B}_{i,\sigma(j)}}) \delta_{n,-m}, \\    \label{dr2}
b_i[n] &= \ep^n b_{\sigma(i)}[n].
\end{align}

We claim that the Poisson algebras ${\mathcal H}'_t(\g)$ and
$\wt{{\mathcal H}}_t(\g)$ are isomorphic. Indeed, for each $i \in \{
1,2,\ldots,\el \}$, choose an element $\ol{i}$ of the orbit of
$\sigma$ in $\Gamma$ corresponding to the simple root $\al_i$. The
long simple roots of $\g$ correspond to simply-transitive orbits of
$\sigma$ in $\Gamma$, while short simple roots correspond to the fixed
points of $\sigma$ in $\Gamma$. Due to formula \eqref{dr2}, we can
choose $b_{\ol{i}}[n], i \in R_l; n \in \Z$, and $\frac{1}{r^\vee}
b_{\ol{i}}[r^\vee n], i \in R_s; n \in \Z$, as the linearly
independent generators of $\wt{{\mathcal H}}_t(\g)$. It is easy to see
that the Poisson brackets between them coincide with the Poisson
brackets between the generators of ${\mathcal H}'_t(\g)$ up to a sign,
which can be removed by an elementary redefinition of the generators.

Now let $^L\G$ be the twisted affine algebra, which is Langlands dual
to $\G$, i.e., the Cartan matrix of $^L\G$ is dual to the Cartan
matrix of $\G$. Consider the quantized enveloping algebra $U_t(^L\G)$
and its Drinfeld's generators $\kappa_{i,n}$ \cite{Dr} (see
Sect.~4.2). In the limit when the level of $U_t(^L\G)$ tends to $0$,
they generate a Poisson algebra ${\mathcal B}_t(^L\G)$. It is easy to
see that ${\mathcal B}_t(^L\G)$ is isomorphic to $\wt{{\mathcal
H}}_t(\g)$, and hence to ${\mathcal H}'_t(\g)$. In view of the last
paragraph of Sect.~4.2, it is natural to make the following
conjecture.

Let $\W'_t(\g)$ be the intersection of $\W_{\ep,t}(\g)$ and ${\mathcal
H}'_t(\g)$. This is a Poisson subalgebra of ${\mathcal H}'_t(\g)$.

\begin{conj}    \label{centertw}
{\em $\W'_t(\g)$ is isomorphic, as a Poisson algebra, to the center of
$U_t(^L\G)$ at the critical level.}
\end{conj}

Some evidence supporting this conjecture will be presented in
Sect.~6.2.

\section{Deformed $\W$--algebras associated to Lie algebras of classical
types} \label{class-types}

\subsection{The fields $\La_i(z)$}

We first define a set $J$ and fields $\La_i(z), i \in J$, for each simple
Lie algebra of classical type.

\subsubsection{The $A_\el$ series}

$J = \{ 1,\ldots,\el+1 \}$.

$$
\La_i(z) = :Y_i(zq^{-i+1} t^{i-1}) Y_{i-1}(zq^{-i} t^i)^{-1}:, \quad
\quad i=1,\ldots,\el+1.
$$

Equivalently,
\begin{align*}
\La_1(z) &= Y_1(z), \\
\La_i(z) &= :\La_i(z) A_{i-1}(zq^{-i+1}t^{i-1})^{-1}:, \quad \quad
i=2,\ldots,\el.
\end{align*}

\subsubsection{The $B_\el$ series}

$J = \{ 1,\ldots,\el,0,\ol{\el},\ldots,\ol{1} \}$.

\begin{align*}
\La_i(z) &= :Y_i(zq^{-2i+2} t^{i-1}) Y_{i-1}(zq^{-2i} t^i)^{-1}:, \quad
\quad i=1,\ldots,\el-1, \\
\La_\el(z) &= :Y_\el(zq^{-2\el+3} t^{\el-1}) Y_\el(zq^{-2\el+1} t^{\el-1})
Y_{\el-1}(zq^{-2\el} t^\el)^{-1}:, \\
\La_0(z) &= \frac{(q+q^{-1})(qt^{-1} - q^{-1} t)}{q^2 t^{-1} - q^{-2} t}
:Y_\el(zq^{-2\el+3} t^{\el-1}) Y_\el(zq^{-2\el-1} t^{\el+1})^{-1}:, \\
\La_{\ol{\el}}(z) &= :Y_{\el-1}(zq^{-2\el+2} t^\el) Y_\el(zq^{-2\el+1}
t^{\el+1})^{-1} Y_\el(zq^{-2\el-1} t^{\el+1})^{-1}:, \\
\La_{\ol{i}}(z) &= :Y_{i-1}(zq^{-4\el+2i+2} t^{2\el-i}) Y_i(zq^{-4\el+2i}
t^{2\el-i+1})^{-1}:, \quad \quad i=1,\ldots,\el-1.
\end{align*}

Equivalently,
\begin{align*}
\La_1(z) &= Y_1(z), \\
\La_i(z) &= :\La_{i-1}(z) A_{i-1}(zq^{-2i+2}t^{i-1})^{-1}:, \quad \quad
i=2,\ldots,\el, \\
\La_0(z) &= \frac{(q+q^{-1})(qt^{-1} - q^{-1} t)}{q^2 t^{-1} - q^{-2} t}
:\La_\el(z) A_\el(zq^{-2\el} t^{\el})^{-1}:, \\
\La_{\ol{\el}}(z) &= \frac{q^2 t^{-1} - q^{-2} t}{(q+q^{-1})(qt^{-1} -
q^{-1} t)} :\La_0(z) A_\el(zq^{-2\el+2} t^{\el})^{-1}:, \\
\La_{\ol{i}}(z) &= :\La_{\ol{i+1}}(z) A_i(zq^{-2(2\el-i-1)}
t^{2\el-i})^{-1}, \quad \quad i=1,\ldots,\el-1.
\end{align*}

\subsubsection{The $C_\el$ series}

$J = \{ 1,\ldots,\el,\ol{\el},\ldots,\ol{1} \}$.

\begin{align*}
\La_i(z) &= :Y_i(zq^{-i+1} t^{i-1}) Y_{i-1}(zq^{-i} t^i)^{-1}:, \quad
\quad i=1,\ldots,\el, \\
\La_{\ol{i}}(z) &= :Y_{i-1}(zq^{-2\el+i-2} t^{2\el-i}) Y_i(zq^{-2\el+i-3}
t^{2\el-i+1})^{-1}:, \quad \quad i=1,\ldots,\el.
\end{align*}

Equivalently,
\begin{align*}
\La_1(z) &= Y_1(z), \\
\La_i(z) &= :\La_{i-1}(z) A_{i-1}(zq^{-i+1}t^{i-1})^{-1}:, \quad \quad
i=2,\ldots,\el, \\
\La_{\ol{\el}}(z) &= :\La_{\el}(z) A_{\el}(zq^{-\el-1}t^{\el})^{-1}:, \\
\La_{\ol{i}}(z) &= :\La_{\ol{i+1}}(z) A_i(zq^{-2\el+i-2} t^{2\el-i})^{-1}:,
\quad \quad i=1,\ldots,\el-1. 
\end{align*}

\subsubsection{The $D_\el$ series}

$J = \{ 1,\ldots,\el,\ol{\el},\ldots,\ol{1} \}$.

\begin{align*}
\La_i(z) &= :Y_i(zq^{-i+1} t^{i-1}) Y_{i-1}(zq^{-i} t^i)^{-1}:, \quad
\quad i=1,\ldots,\el-2, \\
\La_{\el-1}(z) &= :Y_\el(zq^{-\el+2} t^{\el-2}) Y_{\el-1}(zq^{-\el+2}
t^{\el-2}) Y_{\el-2}(zq^{-\el+1} t^{\el-1})^{-1}:, \\
\La_\el(z) &= :Y_\el(zq^{-\el+2} t^{\el-2}) Y_{\el-1}(zq^{-\el}
t^\el)^{-1}:, \\
\La_{\ol{\el}}(z) &= :Y_{\el-1}(zq^{-\el+2} t^{\el-2}) Y_\el(zq^{-\el}
t^\el)^{-1}:, \\
\La_{\ol{\el-1}}(z) &= :Y_{\el-2}(zq^{-\el+1} t^{\el-1})
Y_{\el-1}(zq^{-\el} t^\el)^{-1} Y_\el(zq^{-\el} t^\el):, \\
\La_{\ol{i}}(z) &= :Y_{i-1}(zq^{-2\el+i+2} t^{2\el-i-2}) Y_i(zq^{-2\el+i+1}
t^{2\el-i-1})^{-1}:, \quad \quad i=1,\ldots,\el-2.
\end{align*}

Equivalently,
\begin{align*}
\La_1(z) &= Y_1(z), \\
\La_i(z) &= :\La_{i-1}(z) A_{i-1}(zq^{-i+1}t^{i-1})^{-1}:,
\quad \quad i=2,\ldots,\el, \\
\La_{\ol{\el}}(z) &= :\La_{\el-1}(z) A_{\el}(zq^{-\el+1}t^{\el-1})^{-1}:,
\\
\La_{\ol{\el-1}}(z) &= :\La_{\ol{\el}}(z) A_{\el-1}(zq^{-\el+1}
t^{\el-1})^{-1}: \\
&= :\La_\el(z) A_\el(zq^{-\el+1} t^{\el-1})^{-1}:, \\
\La_{\ol{i}}(z) &= :\La_{\ol{i+1}}(z)
A_i(zq^{-2\el+i+2}t^{2\el-i-2})^{-1}:, \quad \quad 
i=1,\ldots,\el-2.
\end{align*}

\subsection{The first generating field of the deformed $\W$--algebra}

For each Lie algebra of classical type, set
\begin{equation}    \label{t1}
T_1(z) = \sum_{i \in J} \La_i(z).
\end{equation}

In the $A_\el$ case the field $T_1(z)$ coincides with the one obtained in
\cite{SKAO,FF:w,AKOS}.

\begin{thm}    \label{thm1}
{\em For each Lie algebra $\g$ of classical type, the field $T_1(z)$
commutes with the screening operators $S^\pm_i, i=1,\ldots,\el$ and hence
belongs to ${\mathbf W}_{q,t}(\g)$.}
\end{thm}

\begin{proof} As in \cite{FF:w,AKOS}, we need to prove that the
commutator of $T_1(z)$ and $S^\pm_i(w)$, considered as a formal power
series in $z$ and $w$, is a total difference with respect to $w$:
\begin{align*}
[T_1(z),S^+_i(w)] &= {\mathcal D}_{q^2} \cdot R^+_i(z,w), \\
[T_1(z),S^-_i(w)] &= {\mathcal D}_{t^2} \cdot R^-_i(z,w),
\end{align*}
where
$$
{\mathcal D}_{a} \cdot f(w) = \frac{f(w)-f(wa)}{w(1-a)}.
$$
This will immediately imply the statement of the theorem.

This fact has already been proved in \cite{FF:w,AKOS} in the case
$\g=A_\el$ not only for $T_1(z)$, but for all $T_i(z),
i=1,\ldots,\el$. Hence we focus on the remaining series $B_\el,C_\el$,
and $D_\el$.

Note that according to formula \eqref{ay},
$$
[y_i[n],S^\pm_j(w)] =  0, \quad \quad i\neq j.
$$
It is then clear from the explicit formulas for $T_1(z)$ that
nontrivial contribution to the commutator $[T_1(z),S^+_i(w)]$, where
$i=1,\ldots,\el-1$ (resp., $i=1,\ldots,\el-2$) in the case
$\g=B_\el,C_\el$ (resp, $\g=D_\el$) come from the terms
$\La_i(z),\La_{i+1}(z)$, and $\La_{\ol{i}}(z), \La_{\ol{i+1}}(z)$. It
is easy to see from the formula for $T_1(z)$ that the corresponding
contributions coincide with the commutators of $S^\pm_i(w)$ with
$T_1(z)$ and $T_{\el-1}(z)$, respectively, in the case
$\g=A_\el$. Hence they are total differences according to
\cite{FF:w,AKOS}.

For the remaining screening currents: $S^\pm_\el(z)$ in the case
$\g=B_\el,C_\el$, and $S^\pm_{\el-1}(z)$, $S^\pm_\el(z)$ in the case
$\g=D_\el$, it suffices to prove the statement when $\g$ is of types
$B_2$, $C_2$, and $\g=D_3$. In the latter case, the field $T_1(z)$ is
the same as the field $T_2(z)$ for $\g=A_3$, and hence this case
follows again from \cite{FF:w,AKOS}. Thus, we are left with
$\g=B_2,C_2$ and $S^\pm_2(w)$. Here we explain the proof in the most
difficult case $\g=B_2$ and $S^-_2(w)$. In the remaining cases the
statement is proved by a similar computation.

We find using formula \eqref{ay}:
\begin{align*}
\La_2(z) S^-_2(w) &= q^4 \frac{1-\frac{w}{z} t^{-1}}{1-\frac{w}{z} q^4
t^{-1}} :\La_2(z) S^-_2(w):, \quad \quad |z| \gg |w|, \\
S^-_2(w) \La_2(z) &= \frac{1-\frac{z}{w} t}{1-\frac{z}{w} q^{-4}
t} :\La_2(z) S^-_2(w):, \quad \quad |w| \gg |z|.
\end{align*}
This implies the following commutator between the formal power series
$\La_2(z)$ and $ S^-_2(w)$:
\begin{equation}    \label{comm1}
[\La_2(z),S^-_2(w)] =
(q^4-1) \delta \left( \frac{w}{z} q^4 t^{-1} \right) :Y_2(wq^3)
Y_2(wq) Y_1(wt)^{-1} S^-_2(w):.
\end{equation}
We obtain in an analogous way:
\begin{align}    \label{comm2}
[\La_0(z),S^-_2(w)] &= (q^{-2}-q^2) \delta
\left( \frac{w}{z}
q^4 t^{-3} \right) :Y_2(wq^3t^{-2}) Y_2(wq^{-1})^{-1} S^-_2(w): \\
&+ \label{comm3}  (q^2-q^{-2}) \delta
\left( \frac{w}{z} q^2 t^{-1} \right) :Y_2(wq)
Y_2(wq^{-3}t^2)^{-1} S^-_2(w):,
\end{align}
and
\begin{multline}    \label{comm4}
[\La_{\ol{2}}(z),S^-_2(w)] = \\ (q^{-4}-1) \delta \left( \frac{w}{z}
q^2 t^{-3} \right) :Y_1(wt^{-1}) Y_2(wq^{-1})^{-1} Y_2(wq^{-3})^{-1} S^-_2(w):.
\end{multline}
Now recall that
$$
S^-_2(zt) = t^{-2} q^2 :A_2(z) S^-_2(zt^{-1}):
$$
and note that
$$
A_2(z) = :Y_1(z)^{-1} Y_2(zqt^{-1}) Y_2(zq^{-1}t):
$$
(see formula \eqref{expr}). Using these formulas, we find that the
terms \eqref{comm1} and \eqref{comm2} combine into a total
$t^2$--difference, and so do the terms \eqref{comm3} and
\eqref{comm4}. Hence the result is proved.
\end{proof}

\thmref{thm1} proves \conjref{conj1} for all $\g$ of classical types
and their first fundamental representation.

\section{Connection with analytic Bethe Ansatz}

Now we discuss the connection between our formulas for $T_1(z)$ and the
analytic Bethe Ansatz. We expect the same connection to hold not only for
$\g$ of classical types, but for all $\g$.

\subsection{The limit $t \arr 1$}

The specialization of our formula for $T_1(z)$ to $t=1$ coincides with
the Bethe Ansatz formula for the eigenvalues of the transfer-matrix
corresponding to the finite-dimensional representation $V_{\omega_1}$
of $U_q(\G)$ with highest weight $\omega_1$. Let us explain that in
more detail.

For a finite-dimensional representation $V$ of $U_q(\G)$, the Bethe
Ansatz formula for the corresponding eigenvalue $t_V(z)$ is given by a
linear combination of terms of the form $$Y_{i_1}(zq^{a_1}) \ldots
Y_{i_k}(zq^{a_k}),$$ where $Y_i(z)$ are certain functions, such that
the set of weights of these terms is the set of weights of $V$ (see
\cite{Re1,Re2,KS1} for details). If we set $q=1$, then $t_V(z)$
becomes simply the character of $V$ considered as a representation of
$U_q(\g)$. Therefore $t_V(z)$ can be viewed as a ``$q$--character'' of
the representation $V$. These $q$--characters satisfy the following
natural properties: $t_{V \oplus W}(z) = t_V(z) + t_W(z)$, $t_{V
\otimes W}(z) = t_V(z) t_W(z)$, and $t_{V(\lambda)} = t_V(z\lambda)$,
where $V(\la)$ is the standard twist of $V$ by $\lambda \in
\C^\times$.

Systematically, these $q$--characters can be obtained as follows. Each
finite-dimensional representation $V$ of $U_q(\G)$ gives rise to a
generating series $T_V(z)$ of central elements of $U_q(\G)$ at the
critical level \cite{RS} (see also \cite{DE}). In \cite{FR:crit} we
showed that the free field realization of $U_q(A_\el^{(1)})$ gives us
an embedding of the center of $U_q(A_\el^{(1)})$ into the
Heisenberg-Poisson algebra ${\mathcal H}_{q,1}(A_\el)$, and that the
formula for the image of $T_{V_{\omega_i}}(z)$ in ${\mathcal
H}_{q,1}(\g)$ coincides with the formula for $t_{V_{\omega_i}}(z)$
obtained by the Bethe Ansatz. Although the free field realization of
$U_q(\G)$ has not yet been constructed for $\g$ other than $A_\el$, we
expect that it exists for arbitrary $\g$. Furthermore, we expect that
one can reproduce the Bethe Ansatz formulas by applying the free field
realization to the generating series $T_V(z)$ of elements of the
center of $U_q(\G)$. The fact that the limit $t \arr 1$ of our
formulas for $T_1(z)$ agrees with the Bethe Ansatz formula
$t_{V_{\omega_1}}(z)$ therefore provides supporting evidence for
\conjref{center} and a motivation for \conjref{conj1}.

Our derivation of the formula for $T_1(z)$ suggests two features of
the Bethe Ansatz formulas that appear to be new. First, we represent
the terms $\La_i(z)$ in the formula for the eigenvalues as the
products of $Y_1(z)$ corresponding to the highest weight, and the
``step operators'' $A_i(z)^{-1}$ corresponding to the simple
roots. Second, we show that the linear combination $T_1(z)$ of the
terms $\La_i(z)$ is distinguished by the fact that it commutes with
the screening operators $S^+_i$, which are well-defined in the limit
$t \arr 1$.

We expect that \conjref{conj1} has the following generalization: to
each integral dominant weight $\la$ of $\g$ one can attach a field
$T_\la(z)$ from ${\mathbf W}_{q,t}(\g)$, which is the sum of terms
corresponding to weights in the irreducible representation $V_\la$
counted with multiplicity. The limit of $T_\la(z)$ as $t \arr 1$
should coincide with $t_{V_\la}(z)$, and hence one obtains a new
method of singling out the ``$q$--characters'' by their property of
commutativity with the screening operators $S_i^+$. We will discuss
this in more detail in \cite{FR:new}.

In the case of $A_\el$, explicit formulas for other generating fields
$T_i(z), i=2,\ldots\el$, have been found in
\cite{FR:crit,FF:w,AKOS}. Their specialization at $t=1$ coincides with
$q$--characters of the finite-dimensional representations of $A_\el$
with highest weights $\omega_i, i=2,\ldots\el$. In fact, these fields
can be obtained from the field $T_1(z)$ by a ``fusion procedure'',
i.e., by taking the residues in the operator product expansions of the
fields $T_i(z)$. Moreover, for $\g=A_\el$ it seems that all fields
$T_\la(z)$ can be constructed starting from $T_1(z)$ via the fusion
procedure (some examples have been given in \cite{FR:dca}), and so
${\mathbf W}_{q,t}(A_\el)$ is the smallest subalgebra of the DCA
$(\hh,\pi_0)$, containing $T_1(z)$. For general $\g$, it would be
interesting to analyze which fields $T_\la(z)$ can be obtained from
$T_1(z)$ by the fusion procedure.

We also remark that our formulas above and formulas from
\cite{FF:w,AKOS} suggest that the fields $T_\la(z)$ should have the
following symmetry: if we replace $q \arr q^{-1}, t \arr t^{-1}$, and
$Y_i(z) \arr Y_i(z)^{-1}$, then $T_\la(z)$ becomes $T_{\ol{\la}}(z)$,
where $\ol{\la}$ is the highest weight of the representation dual to
$V_\la$, up to an overall multiplication of $z$ by a factor $q^a
t^b$. In our formulas above, this factor equals $q^{r^\vee h^\vee}
t^{-h}$.

\subsection{The limit $q \arr \ep$}

For simply-laced $\g$, we have the duality $\W_{q,t}(\g) \simeq
\W_{t,q}(\g)$. Hence $\W_{1,t}(\g)$ is isomorphic to
$\W_{t,1}(\g)$. Thus, without loss of generality in this subsection we
can focus on nonsimply-laced $\g$.

Recall that according to \conjref{centertw}, the subalgebra
$\W'_t(\g)$ of $\W_{\ep,t}(\g)$ is isomorphic to the center of
$U_t(\GL)$. Just as in the non-twisted case (see the previous
subsection), to every finite-dimensional representation $V$ of
$U_t(\GL)$ corresponds a Bethe Ansatz formula -- its
``$t$--character''. On the other hand, it also gives rise to a
generating series $T_V(z)$ of central elements of $U_t(\GL)$. We
expect that $U_t(\GL)$ has a free field realization, which embeds the
center of $U_t(\GL)$ into a Heisenberg-Poisson algebra, so that the
image of $T_V(z)$ coincides with the formula for the $t$--character of
$V$.

Therefore we expect that the formulas for the elements of
$\W_{q,t}(\g)$ that lie in $\W'_t(\g)$ in the limit $q \arr 1$,
coincide with the corresponding Bethe Ansatz formulas. Examples of the
latter formulas are known in the literature \cite{Re2,KS2}, and we can
compare them with our formulas.

Let us first consider the formula for $T_1(z)$ in the case of
$B_\el$. In this case, all simple roots, except for the $\el$th root,
are long and the $\el$th root is short. Hence an element of
$\W_{q,t}(B_\el)$ lands in $\W'_t(B_\el)$ as $q \arr i$ if $Y_\el(z)$
appears in it only through combination $Y_\el(z) Y_\el(zq^2)$, giving
$Y_\el(z) Y_\el(-z)$ after the specialization $q=i$. We see from our
formula that this is so for $T_1(z)$. Furthermore, the term $\La_0(z)$
vanishes at $q=i$, so we are left with $2\el$ terms. Comparison with
\cite{Re2,KS2} shows perfect agreement between the resulting formula
for $T_1(z)$ and the Bethe Ansatz formula for the $2\el$--dimensional
representation of $U_t(A^{(2)}_{2\el-1})$ (note that
$^L\!(B^{(1)}_\el) = A^{(2)}_{2\el-1}$). This provides supporting
evidence for \conjref{centertw}.

Thus, $T_1(z)$ ``interpolates'' between the $q$--character of the
$(2\el+1)$--dimensional representation of $U_q(B^{(1)}_\el)$ (in the
limit $t \arr 1$) and the $t$--character of the $2\el$--dimensional
representation of $U_t(A^{(2)}_{2\el-1})$.

In the case of $C_\el$, the short simple roots are
$\al_1,\ldots,\al_{\el-1}$, so for an element of $\W_{q,t}(C_\el)$ to
land in $\W'_t(C_\el)$ as $q \arr i$, the fields $Y_i(z)$ should
appear in it only in combination $Y_i(z) Y_i(zq^2)$ for
$i=1,\ldots,\el-1$. Inspection of our formula for $T_1(z)$ shows that
it does not satisfy this condition. However, the specializations to
$q=i$ of other fields from ${\mathbf W}_{q,t}(C_\el)$ may well satisfy
it, in which case they should coincide with the $t$--characters of
certain representations of $U_t(D^{(2)}_{\el+1})$, since
$^L\!(C^{(1)}_\el) = D^{(2)}_{\el+1}$.

We have constructed such a field in the case $\g=C_2$. In the notation
of the previous subsection, this is the field $T_{2\omega_1}(z)$. When
$t=1$ it coincides with the $q$--character of the representation
$V_{2\omega_1}$ of $U_q(C_2^{(1)})$, while at $q=i$ it coincides with
the $t$--character of the representation $V_{\omega_1}$ of
$U_t(D_3^{(2)})$. We expect that the same is true when $\g=C_\el$.

It would be interesting to see how this duality works for general
representations, and whether it actually sets up a correspondence
between certain finite-dimensional representations of $U_q(\G)$ and
$U_t(\GL)$.

\subsection{The limit $q \arr 1$}

Finally, let us consider the special case $q=1$ of our formula for
$T_1(z)$. We claim that it is closely related to the Bethe Ansatz
formula for the first fundamental representation of $U_t(\G^\vee)$,
where $\G^\vee = ^L\!\!(\widehat{^L\!\g})$ (e.g., $(B_\el^{(1)})^\vee
= D_{\el+1}^{(2)}$, and $(C_\el^{(1)})^\vee = A_{2\el-1}^{(2)}$). The
latter can be found in \cite{Re2,KS2}. The only difference is that in
the models associated to twisted affine algebras, Bethe Ansatz
formulas contain terms of the form $Y_i(\ep^k z t^a)$. Our formulas
coincide with formulas for the eigenvalues modulo the $\ep$ factors,
which are absent in our formulas when $t=1$ (these factors do appear
in the limit $q \arr \ep$, see the previous subsection). This
indicates that our formulas at $q=1$ are $t$--characters of
finite-dimensional representations of some algebra closely related to
$U_t(\G^\vee)$.

It is again interesting to analize explicitly how finite-dimensional
representations of the dual pairs of affine algebras, $\G$ and
$\G^\vee$, get connected in the deformed $\W$--algebra
$\W_{q,t}(\g)$. This is clear for simply-laced $\g$, but one can see
rather intriguing effects for nonsimply-laced $\g$, similar to those
considered in the previous subsection.

For instance, in the case $\g=B_\el$ the deformed $\W$--algebras
apparently interpolates between finite-dimensional representations of
$B_\el^{(1)}$ and $D_{\el+1}^{(2)}$. The first fundamental
representation of $B_\el^{(1)}$ has dimension $2\el+1$, while the
first fundamental representation of $D_{\el+1}^{(2)}$ has dimension
$2\el+2$. We see from the formula for $T_1(z)$ above that at $t=1$
there are indeed $2\el+1$ terms corresponding to the weight spaces in
the fundamental representation of $B_\el^{(1)}$. On the other
hand, at $q=1$, one of the terms gets doubled, and there appear
$2\el+2$ terms, corresponding to the weight spaces in the fundamental
representation of $D_{\el+1}^{(2)}$.

Similar effect probably occurs for the dual pairs
$(F_4^{(1)},E_6^{(2)})$ ($26$-- and $27$--dimensional fundamental
representations, respectively) and $(G_2^{(1)},D_4^{(3)})$ ($7$--
and $8$--dimensional fundamental representations, respectively).

In the case $\g=C_\el$ our formula for $T_1(z)$ connects the
$2\el$--dimensional fundamental representation of $C_\el^{(1)}$ and
the $2\el$--dimensional fundamental representation of
$A_{2\el-1}^{(2)}$.

\section{The case of $A_{2\el}^{(2)}$}

Up to now, we have not discussed the series $A_{2\el}^{(2)}$ of
self-dual twisted affine algebras. In this section we define the
deformed $\W$--algebra associated to $A_{2\el}^{(2)}$ and its free
field realization. One may also view it as the $\W$--algebra
associated to the non-reduced root system $BC_\el$. The results of
this section generalize the results of Brazhnikov and Lukyanov
\cite{BL}, who constructed the deformed $\W$--algebra associated to
$A_2^{(2)}$.

First we define the $\el \times \el$ symmetric matrix $\B(q,t)$ associated
to $A_{2\el}^{(2)}$ as follows: 
\begin{align*}
\B_{ij}(q,t) &= [2]_q \B_{ij}^A(q^2,t), \quad \quad (i,j) \neq (\el,\el), \\
\B_{\el\el}(q,t) & = [2]_q (q^2t^{-1} - 1 + q^{-2}t),
\end{align*}
where $(\B^A_{ij}(q,t))$ is the $\B$--matrix associated to $A_\el$.

By definition, the Heisenberg algebra $\HH(\A)$ has generators
$a_i[n], i=1,\ldots,\el;$ $n \in \Z$, and relations \eqref{a} with
$B(q,t)$ as above. We define the ``dual'' generators $y_i[n]$ by
formula \eqref{y} with $r_i=2, \forall i$, and the operators $e^{Q_j}$
by formula \eqref{qj}, where $B_{ij}=B_{ij}(1,1)$.

Now set 
$$
A_i(z) = t^2 q^{-4 + 2a_i[0]} :\exp \left( \sum_{m\neq 0} a_i[m] z^{-m}
\right):,
$$
$$
Y_i(z) = t^{i(2\el+1-i)} q^{-2i(2\el+1-i) + 2 y_i[0]} :\exp \left(
\sum_{m\in\Z} y_i[m] z^{-m} \right):.
$$

The screening operators are defined by formulas \eqref{s+} and \eqref{s-},
in which we set $\rr_i=2, \forall i$.

Now we define the DCA ${\mathbf W}_{q,t}(\A)$ as the maximal
subalgebra of $({\mathbf H}_{q,t}(\A),\pi_0)$, which commutes with the
operators $S^-_i, i=1,\ldots,\el$. We define the deformed
$\W$--algebra $\W_{q,t}(\A)$ as the associative algebra, topologically
generated by the Fourier coefficients of fields from ${\mathbf
W}_{q,t}(\A)$. We expect an analogue of \conjref{conj1} to hold in the
$\A$ case.

Here is the explicit formula for the field $T_1(z)$:
$$
T_1(z) = \sum_{i \in J} \La_i(z),
$$
where $J = \{ 1,\ldots,\el,0,\ol{\el},\ldots,\ol{1} \}$ and

\begin{align*}
\La_i(z) &= :Y_i(zq^{-2(i+1)} t^{i-1}) Y_{i-1}(zq^{-2i} t^i)^{-1}:, \quad
\quad i=1,\ldots,\el, \\
\La_0(z) &= \frac{1-q^{-2}t^{-1}}{1-q^2t^{-1}} :Y_\el(zq^{-2\el} t^\el)
Y_\el(zq^{-2(\el+1)} t^{\el+1})^{-1}:, \\
\La_{\ol{i}}(z) &= q^{-4} :Y_{i-1}(zq^{-2(2\el-i+1)} t^{2\el-i+1})
Y_i(zq^{-2(2\el-i+2)} t^{2\el-i+2})^{-1}:, \quad \quad i=1,\ldots,\el-1.
\end{align*}

Equivalently,
\begin{align*}
\La_1(z) &= Y_1(z), \\
\La_i(z) &= :\La_{i-1}(z) A_{i-1}(zq^{-2(i-1)}t^{i-1})^{-1}:, \quad \quad
i=2,\ldots,\el, \\
\La_0(z) &= \frac{1-q^{-2}t^{-1}}{1-q^2t^{-1}} :\La_\el(z) A_\el(zq^{-2\el}
t^\el)^{-1}:, \\
\La_{\ol{\el}}(z) &= \frac{1-q^{-2}t}{1-q^2t} :\La_0(z)
A_\el(zq^{-2(\el+1)} t^{\el+1})^{-1}:, \\
\La_{\ol{i}}(z) &= :\La_{\ol{i+1}}(z) A_i(zq^{-2(2\el-i+1)}
t^{2\el-i+1})^{-1}, \quad \quad i=1,\ldots,\el-1.
\end{align*}

Again, we find a perfect agreement between the specialization of this
formula to $q=i$ and the formula obtained by analytic Bethe Ansatz in
the $A_{2\el}^{(2)}$ integrable model (see \cite{Re2,KS2})\footnote{in
the case of $A^{(2)}_2$ this has been observed in \cite{BL}}, in
agreement with our general scheme outlined above.

Note however, that in contrast with the general case, $T_1(z)$ (and
probably other fields as well) does not commute with the second set of
screening operators, $S^+_i$.

\section{The scaling limit}

\subsection{The scaling limit of the exchange relations and factorized
scattering}

In this subsection we show that in the scaling limit the function
$S_{ij}(x) = S_{Y_i,Y_j}(x)$ given by formula \eqref{theta} becomes
the $S$--matrix of the $i$th and $j$th particles of the affine Toda
field theory associated to $\G$ \cite{BCDS,DGZ,CDS}, in the integral
form given by T.~Oota \cite{O}.

Let us set
$$
z=e^{-i\theta_1 \ep}, \quad w=e^{-i\theta_2 \ep},
$$
$$
q=\exp \frac{B}{2r^\vee h^\vee} \pi \ep, \quad \quad t=\exp \frac{B-2}{2h}
\pi \ep.
$$
Put $\theta=\theta_1-\theta_2$. Then in the limit $\ep \arr 0$ with
$\theta_i$ and $B$, the function $S_{ij}(x)$ becomes:
\begin{multline*}
S_{ij}(\theta) = \\ \exp \left( - 4 \int_{-\infty}^\infty \frac{dk}{k}
e^{ik\theta} \cdot \on{sinh} \frac{2-B}{2h} \pi k \cdot \on{sinh}
\frac{B}{2r^\vee h^\vee} \pi k \cdot M_{ij} \left( e^{\frac{B}{2r^\vee
h^\vee} \pi k},e^{\frac{B-2}{2h} \pi k} \right) \right),
\end{multline*}
(the integral should be understood as the principal value).

This formula coincides with Oota's integral formula for the
$S$--matrix of the affine Toda field theory \cite{O}.\footnote{Note
that a factor of $i$ is missing in front of the integral in \cite{O},
and that the matrix $(C_{ab})$ used in \cite{O} is the transpose of
the Cartan matrix according to the conventions of \cite{Kac}. We thank
T.~Oota for clarifying this point to us.} This means that the scaling
limit of our deformed $\W$--algebra $\W_{q,t}(\g)$ can be viewed as
the Faddeev-Zamolodchikov algebra of the Toda theory. This has been
suggested by Lukyanov in the $A_\el$ case \cite{Lu2,Lu3} (see also
\cite{Lu1}). The scaling limit of $\W_{q,t}(A_\el)$ has also been
studied in \cite{Hou}. The scaling limit of our free field realization
can be used to obtain explicit formulas for form-factors in general
affine Toda field theories along the lines of \cite{Lu2,Lu3,BL}.

On the other hand, the functions $S_{ij}(\theta)$ describe (at least
in the simply-laced case) the scattering amplitudes for the lowest
breather particles in the affine Toda field theories with imaginary
coupling constant (these particles are bound states of solitons
analogous to the breather particles of the sine-Gordon model). This
fact is connected with the particle-breather duality of Toda theory
\cite{G,J}. Note that the Toda theories with imaginary coupling
constant are non-unitary (except for $\g=A_1$). In order to have a
physically meaningful quantum field theory, one has to restrict those
theories to RSOS scattering states.

The functions \eqref{theta} are elliptic generalizations of the
$S$--matrices of the affine Toda field theories. Hence they can also
be interpreted as the $S$--matrices of the RSOS model, whose scaling
limit gives the corresponding restricted Toda theory with imaginary
coupling constant. The algebra $W_{q,t}(\g)$ can therefore be
interpreted as the Faddeev-Zamolodchikov algebra for this RSOS model
(in the case $\g=A_1$ this has been suggested in \cite{Lu2}). The
embedding of $\W_{q,t}(\g)$ into ${\mathcal H}_{q,t}(\g)$ give us a
bosonization of this Faddeev-Zamolodchikov algebra.

\subsection{The scaling limit of the screening currents}

Set
\begin{equation}    \label{sclim}
x=e^{2\pi i \tau u}, \quad \quad q=e^{2\pi i \tau}, \quad \quad
t=e^{2\pi i \tau \beta}.
\end{equation}

Denote by $G_{ij}^\pm(w/z)$ be the function appearing in the quadratic
relation between $S_i^\pm(z)$ and $S_j^\pm(w)$ from Sect.~3.5:
$$
S^\pm_i(z) S^\pm_j(w) = G^\pm_{ij}\left( \frac{w}{z} \right) S^\pm_j(w)
S^\pm_i(z).
$$

Consider the limit of $G^\pm_{ij}(x)$ as $\tau \arr + i0$ with $x, q$,
and $t$ given by \eqref{sclim}. Using the modular properties of the
function $\theta(z;a)$, we find:
$$
\lim_{\tau \arr + i0} \frac{\theta(e^{2\pi i
\tau(u+\al)};e^{2\pi i \tau})}{\theta(e^{2\pi i
\tau(u+\gamma)};e^{2\pi i \tau})} = \frac{\on{sin}
\pi(u+\al)}{\on{sin} \pi(u+\gamma)}.
$$
Applying this formula, we obtain the following limits for the functions
$G^-_{ij}(x)$:
\begin{align}
G^-_{kk}(x) &\arr \frac{\sin \frac{\pi}{2\beta}(u-2)}{\sin
\frac{\pi}{2\beta}(u+2)}, \\
G^-_{km}(x) &\arr \frac{\sin \frac{\pi}{2\beta}(u+\beta-B_{km})}{\sin
\frac{\pi}{2\beta}(u+\beta+B_{km})}, \quad \quad k\neq m.
\end{align}

By shifting the variables $u_k, k=1,\ldots,\el$, in the obvious way,
we obtain the following relations between the scaling limits $s^-_k(u)$
of the screening currents:
\begin{equation}    \label{s-1}
s^-_k(u) s^-_m(v) = \frac{\sin \frac{\pi}{2\beta}(v-u-B_{km})}{\sin
\frac{\pi}{2\beta}(v-u+B_{km})} s^-_m(v) s^-_k(u).
\end{equation}
If we set $z=e^{\pi iu/\beta}, w=e^{\pi iv/\beta}$, and $q=e^{\pi i/\beta}$, then
the function $g^-_{km}$ appearing in the right hand side of formula
\eqref{s-} becomes
$$
g^-_{km} = \frac{zq^{B_{ij}} - w}{z - w q^{B_{ij}}},
$$
and we see that relations \eqref{s-1} coincide with the Drinfeld
relations \cite{Dr} in the subalgebra $U_q(\widehat{\n})$ of
$U_q(\G)$. 

Analogously, we obtain the following limits for the functions
$G^+_{km}(x)$:
\begin{align}
G^+_{kk}(x) &\arr \frac{\sin \frac{\pi}{2r_{kk}}(u-2\beta)}{\sin
\frac{\pi}{2r_{kk}}(u+2\beta)}, \\
G^+_{km}(x) &\arr \frac{\sin \frac{\pi}{2r_{km}}(u+r_{km}+\beta)}{\sin
\frac{\pi}{2r_{km}}(u+r_{km}-\beta)}, \quad \quad k\neq m, B_{km} \neq
0.
\end{align}
By shifting the variables $u_k, k=1,\ldots,\el$, we obtain the
following relations between the scaling limits $s^+_k(u)$
of the screening currents:
\begin{align}    \label{s+1}
s^+_k(u) s^+_k(v) &= \frac{\sin \frac{\pi}{2r_{ii}}(v-u-2\beta)}{\sin
\frac{\pi}{2r_{ii}}(v-u+2\beta)} s^+_k(v) s^+_k(u), \\
\label{s+2}    s^+_k(u) s^+_m(v)
&= \frac{\sin \frac{\pi}{r_{km}}(v-u+\beta)}{\sin
\frac{\pi}{r_{km}}(v-u-\beta)} s^+_m(v) s^+_k(u), \quad \quad
k\neq m, B_{km} \neq 0.
\end{align}
Denote by $g^+_{ij}$ the function appearing in the right hand sides of
formulas \eqref{s+1}, \eqref{s+2}. Set $z=e^{\pi iu/r^\vee},
w=e^{\pi iv/r^\vee}$, $t=e^{-\pi i\beta/r^\vee}$. Then the functions
$g^+_{km}$ become
\begin{align*}
g^+_{kk} &= \frac{zt^2 - w}{z - w t^2}, \quad \quad k \in R_l, \\
g^+_{km} &= \frac{zt^{-1} - w}{z - w t^{-1}}, \quad \quad k \neq m \in
R_l, B_{km} \neq 0, \\
g^+_{kk} &= \frac{z^{r^\vee} t^{2r^\vee} - w^{r^\vee}}{z^{r^\vee} -
w^{2r^\vee} t^{2r^\vee}}, \quad \quad k \in R_s, \\
g^+_{km} &= \frac{z^{r^\vee}t^{-r^\vee} - w^{r^\vee}}{z^{r^\vee} -
w^{r^\vee} t^{-r^\vee}}, \quad \quad k \neq m, k \in R_s, B_{km}
\neq 0,
\end{align*}
where $R_l$ and $R_s$ denote the sets of long and short simple
roots, respectively. It is easy to recognize in these formulas
the Drinfeld relations in the subalgebra $U_t(\widehat{\n})$ of
$U_t(\GL)$.

Thus, the relations of Sect.~3.5 between the screening currents
$S^-_i(z)$ (resp., $S^+_i(z)$) can be considered as elliptic analogues of
Drinfeld relations in $U_q(\G)$ (resp., $U_t(\GL)$).

\section{Vertex operators}

Vertex operators are constructed from the fields $Y_i(z)$ in the same
way as the screening currents are constructed from the fields
$A_i(z)$.

Introduce the operators $e^{Q_{\om_i}}, i=1,\ldots,\el$, acting from
$\pi_\mu$ to $\pi_{\mu+\beta\om_i}$, which satisfy commutation relations
$$
[a_i[n],e^{Q_{\om_j}}] = \rr_i \delta_{ij} \beta \delta_{n,0} e^{Q_{\om_j}}.
$$

Let
\begin{equation}    \label{vm+}
v^+_i[m] = \frac{y_i[m]}{q^{m \rr_i}-q^{-m \rr_i}}, \quad m\neq 0, \quad
\quad v^+_i[0] = y_i[0]/\rr_i,
\end{equation}
\begin{equation}    \label{vm-}
v^-_i[m] = \frac{y_i[m]}{t^m-t^{-m}}, \quad m\neq 0,
\quad \quad v^-_i[0] = y_i[0]/\beta.
\end{equation}

Now define the {\em fundamental vertex operators} by the formulas
\begin{align}    \label{v+}
V_i^+(z) & = e^{Q_{\om_i}/\rr_i} z^{v^+_i[0]} :\exp \left( - \sum_{m\neq 0}
v^+_i[m] z^{-m} \right):,\\    \label{v-}
V_i^-(z) & = e^{-Q_{\om_i}/\beta} z^{-v^-_i[0]} :\exp \left( \sum_{m\neq 0}
v^-_i[m] z^{-m} \right):.
\end{align}

They satisfy the difference equations:
\begin{equation}    \label{v1}
V^+_i(zq^{\rr_i}) = t^{-2(\rho^\vee,\om_i)} q^{2r^\vee (\rho,\om_i)}
:Y_i(z) V^+_i(zq^{-\rr_i}):,
\end{equation}
and
\begin{equation}    \label{v2}
V^-_i(zt^{-1}) = t^{-2(\rho^\vee,\om_i)} q^{2r^\vee (\rho,\om_i)} :Y_i(z)
V^-_i(zt):.
\end{equation}

For $\g=A_1$ these operators have been constructed in \cite{LP1,LP2}, and
for $\g=A_\el$ they have been constructed in \cite{AJMP}.

In the conformal limit, $q \arr 1, t=q^\beta$, the fields $V^-_i(z)$ and
$V^+_i(z)$ become the bosonizations of the highest weight components of
primary fields of $\W_\beta(\g)$. These primary fields correspond to the
fundamental representations of $\g$ (generalizations of $\Phi_{1,2}(z)$)
and fundamental representations of $^L\!\g$ (generalizations of
$\Phi_{2,1}(z)$), respectively. Other components of these primary fields,
corresponding to other weight components of the fundamental
representations, are obtained by applying to them the screening operators
of the same type ($-$ or $+$). Langlands duality interchanges the $-$ and
$+$ vertex operators.

In the important works \cite{LP1,LP2,AJMP} it was shown that similar
picture holds in the deformed situation in the case $\g=A_\el$. In
those papers a complete bosonization of the vertex operators of the
SOS models associated to $A^{(1)}_\el$ and their restricted
counterparts had been obtained. The highest weight components of those
vertex operators are $V^-_i(z)$ (in our notation), and the rest are
obtained by applying to them the deformed screening operators
$S^-_j$. These are the ``type I operators''. It was explained in
\cite{AJMP} that these operators are bosonizations of the
half-infinite transfer-matrices on the lattice for the $A^{(1)}_\el$
face model. We expect that the fields $V^-_i(z)$ associated to general
Lie algebras can be used in the same fashion to obtain the
bosonization of type I operators in the corresponding SOS and RSOS
models \cite{JMO}.

The ``type II operators'' should be obtained in the same way, starting from
$V^+_i(z)$. However, note that in the conformal limit the components of the
type II operators correspond to weight spaces in the fundamental
representations of $^L\!\g$, rather than $\g$.

\begin{rem} Note that in the classical limits $t \arr 1$ (respectively, $q
\arr \ep$), the vertex operators $V^+_i(z)$ (respectively, $V^-(z)$) become
Baxter's $Q$--operators of the corresponding integrable model.\qed
\end{rem}

\section*{Appendix A. Poisson algebras $\W_{q,1}(\g)$}

In this section we make explicit computations in the Poisson algebra
$\W_{q,1}(\g)$, when $\g$ is a simple Lie algebra of classical type.

Recall that $M(q,t) = D(q,t) \B(q,t)^{-1} D(q,t)$. For classical
Lie algebras, these matrices are given in Appendix C. Set
\begin{align*}
{\mathcal B}_{ij}(x) &= \sum_{m\in\Z} (q^m-q^{-m}) \B_{ij}(q^m,1) x^m, \\
{\mathcal M}_{ij}(x) &= \sum_{m\in\Z} (q^m-q^{-m}) M_{ij}(q^m,1) x^m.
\end{align*}

We have the following Poisson brackets:
\begin{align*}
\{ A_i(z),A_j(w) \} &= {\mathcal B}_{ij} \left( \frac{w}{z} \right) A_i(z)
A_j(w), \\
\{ Y_i(z),Y_j(w) \} &= {\mathcal M}_{ij} \left( \frac{w}{z} \right) Y_i(z)
Y_j(w).
\end{align*}

\bigskip

\noindent A.1. {\bf $A_\el$ case.}
In this section we recall the results of \cite{FR:crit}.

We have:
$$\La_i(z) = Y_i(zq^{-i+1}) Y_{i-1}(zq^{-i})^{-1}, \quad \quad
i=1,\ldots,\el+1.$$

The generators of $\W_{q,1}(A_\el)$ are
$$\s_i(z) = \sum_{1\leq j_1 < \ldots < j_i\leq n+1} \La_{j_1}(z)
\La_{j_2}(zq^2) \ldots \La_{j_{i-1}}(zq^{2(i-2)}) \La_{j_i}(zq^{2(i-1)}),
\quad i=1,\ldots,\el+1.$$ We have: $\s_{\el+1}(z) = 1$.

The Poisson brackets are \cite{FR:crit}:
\begin{align*}
\{ \s_i(z),\s_j(w) \} &= {\mathcal M}_{ij}\left( \frac{wq^{j-i}}{z}
\right) \s_i(z) \s_j(w) \\ &+ \sum_{p=1}^i \delta\left(
\frac{w}{zq^{2p}} \right) \s_{i-p}(w) \s_{j+p}(z)
\\ &- \sum_{p=1}^i \delta\left( \frac{wq^{2(j-i+p)}}{z} \right)
\s_{i-p}(z) \s_{j+p}(w),
\end{align*}
if $i\leq j$ and $i+j\leq \el+1$; and
\begin{align*}
\{ \s_i(z),\s_j(w) \} &= {\mathcal M}_{ij}\left(
\frac{wq^{j-i}}{z} \right) \s_i(z) \s_j(w) \\  &+ \sum_{p=1}^{N-j}
\delta\left( \frac{w}{zq^{2p}} \right) \s_{i-p}(w) \s_{j+p}(z)
\\ &- \sum_{p=1}^{N-j} \delta\left( \frac{wq^{2(j-i+p)}}{z} \right)
\s_{i-p}(z) \s_{j+p}(w),
\end{align*}
if $i\leq j$ and $i+j>\el+1$.

\bigskip

\noindent A.2. {\bf $B_\el$ case.}
This and the next two subsections are taken from \cite{FR:BC}.

We have:
\begin{align*}
\La_i(z) &= Y_i(zq^{-2i+2}) Y_{i-1}(zq^{-2i})^{-1}, \quad \quad
i=1,\ldots,\el-1, \\ \La_\el(z) &= Y_\el(zq^{-2\el+3}) Y_\el(zq^{-2\el+1})
Y_{\el-1}(zq^{-2\el} t^\el)^{-1}, \\ \La_0(z) &= Y_\el(zq^{-2\el+3})
Y_\el(zq^{-2\el-1})^{_1}, \\ \La_{\ol{\el}}(z) &= Y_{\el-1}(zq^{-2\el+2})
Y_\el(zq^{-2\el+1})^{-1} Y_\el(zq^{-2\el-1})^{-1}, \\ \La_{\ol{i}}(z) &=
Y_{i-1}(zq^{-4\el+2i+2}) Y_i(zq^{-4\el+2i})^{-1}, \quad \quad
i=1,\ldots,\el-1.
\end{align*}

Let
$$\s_i(z) = \sum_{\{ j_1,\ldots,j_i \} \in S} \La_{j_1}(z) \La_{j_2}(zq^4)
\ldots \La_{j_i}(zq^{4i-4}), \quad \quad i=1,\ldots,\el-1,$$ where $S$ is the
set of $\{ j_1,\ldots,j_i \}$, such that $j_\al<j_{\al+1}$ or
$j_\al=j_{\al+1}=\el+1, \al=1,\ldots,i-1$.

The formula for $\s_\el(z)$ can be found in \cite{FR:crit}.

\begin{rem} In the formulas above, $q$ corresponds to $q^{1/2}$ of
\cite{FR:crit}.\qed
\end{rem}
\medskip

The Poisson brackets are:

\begin{align*}
\{ \La_i(z),\La_i(w) \} &= {\mathcal M}_{11}\left( \frac{w}{z} \right)
\La_i(z) \La_i(w), \\
\{ \La_i(z),\La_j(w) \} &= {\mathcal M}_{11}\left( \frac{w}{z} \right)
\La_i(z) \La_j(w) +
\left( \delta \left( \frac{w}{zq^4} \right) - \delta \left( \frac{w}{z}
\right) \right. \\ &+ \left. \delta \left( \frac{w}{zq^{4\el-4i+2}}
\right) \delta_{i,\ol{j}} - \delta \left( \frac{w}{zq^{4\el-4i-2}}
\right) \delta_{i,\ol{j}} \right) \La_i(z) \La_j(w),
\end{align*}
if $i<j$.

These formulas imply for $\el>2$:
\begin{align*}
\{ T_1(z),T_1(w) \} &= {\mathcal M}_{11} \left( \frac{w}{z} \right) T_1(z)
T_1(w) \\ &+ \delta \left( \frac{w}{zq^4} \right) T_2(z) - \delta \left(
\frac{wq^4}{z} \right) T_2(w) \\ &+ \delta \left( \frac{w}{zq^{4\el-2}}
\right) - \delta \left( \frac{w^{4\el-2}}{z} \right).
\end{align*}
In the case of $B_2$ the formulas for the Poisson brackets can be read
off the formulas for the Poisson brackets in the case of $C_2$ given in
Sect.~A.4.

\bigskip

\noindent A.3. {\bf $C_\el$ case.}
We have:
\begin{align*}
\La_i(z) &= Y_i(zq^{-i+1}) Y_{i-1}(zq^{-i})^{-1}, \quad \quad
i=1,\ldots,\el, \\
\La_{\ol{i}}(z) &= Y_{i-1}(zq^{-2\el+i-2}) Y_i(zq^{-2\el+i-3})^{-1},
\quad \quad i=1,\ldots,\el.
\end{align*}

Let
$$\s_i(z) = \sum_{\{ j_1,\ldots,j_i \} \in S} \La_{j_1}(z) \La_{j_2}(zq)
\ldots \La_{j_i}(zq^{i-1}), \quad \quad i=1,\ldots,\el,$$ where $S$ is the
set of $\{ j_1,\ldots,j_i \}$, such that $1\leq j_1 < \ldots < j_i\leq 2\el$
and if $j_\al = k, j_\beta = 2\el+1-k$ for some $k=1,\ldots,\el$, then $l\leq
\el+\al-\beta$.

\begin{rem} In the formulas above, $q$ corresponds to $q^{1/2}$ of
\cite{FR:crit}.\qed
\end{rem}
\medskip

The Poisson brackets are:

\begin{align*}
\{ \La_i(z),\La_i(w) \} &= {\mathcal M}_{11}\left( \frac{w}{z} \right) \La_i(z)
\La_i(w), \\
\{ \La_i(z),\La_j(w) \} &= {\mathcal M}_{11}\left( \frac{w}{z} \right) \La_i(z)
\La_j(w) +
\left( \delta \left( \frac{w}{zq^2} \right) - \delta \left( \frac{w}{z}
\right) \right. \\ &+ \left. \delta \left( \frac{w}{zq^{2\el-2i+4}} \right)
\delta_{i,\ol{j}} - \delta \left( \frac{w}{zq^{2\el-2i+2}} \right)
\delta_{i,\ol{j}} \right) \La_i(z) \La_j(w),
\end{align*}
if $i<j$.

These formulas imply:
\begin{align*}
\{ T_1(z),T_1(w) \} &= {\mathcal M}_{11} \left( \frac{w}{z} \right) T_1(z)
T_1(w) \\ &+ \delta \left( \frac{w}{zq^2} \right) T_2(z) - \delta \left(
\frac{wq^2}{z} \right) T_2(w) \\ &+ \delta \left( \frac{w}{zq^{2\el+2}}
\right) - \delta \left( \frac{w^{2\el+2}}{z} \right).
\end{align*}

\bigskip

\noindent A.4. {\bf Poisson brackets for $C_2$.}
In this subsection we give a compute description of the Poisson algebra
$\W_{q,1}(C_2)$.

\begin{align*}
\{ T_1(z),T_1(w) \} &= \sum_{m\in\Z}
\frac{(q^m-q^{-m})(q^{2m}+q^{-2m})}{q^{3m}+q^{-3m}} \left( \frac{w}{z}
\right)^m T_1(z) T_1(w) \\ &+ \delta \left( \frac{w}{zq^2} \right) T_2(z)
-\delta \left( \frac{wq^2}{z} \right) T_2(w) \\ &+ \delta \left(
\frac{w}{zq^6} \right) - \delta \left( \frac{wq^6}{z} \right). \\
\{ T_1(z),T_2(w) \} &= \sum_{m\in\Z} \frac{q^{2m}-q^{-2m}}{q^{3m}+q^{-3m}}
\left( \frac{wq}{z} \right)^m T_1(z) T_2(w) \\ &+ \delta \left(
\frac{w}{zq^4} \right) T_1(w) -\delta \left( \frac{wq^6}{z} \right)
T_1(wq^2), \\
\{ T_2(z),T_2(w) \} &= \sum_{m\in\Z}
\frac{(q^{2m}-q^{-2m})(q^m+q^{-m})}{q^{3m}+q^{-3m}} \left( \frac{w}{z}
\right)^m T_2(z) T_2(w) \\ &+ \delta \left( \frac{w}{zq^4} \right) \left(
T_1(z) T_1(zq^2) - T_2(zq^2) \right) \\ &- \delta \left( \frac{wq^4}{z}
\right) \left( T_1(w) T_1(wq^2) - T_2(wq^2) \right) \\ &+ \delta \left(
\frac{w}{zq^6} \right) - \delta \left( \frac{wq^6}{z} \right).
\end{align*}

\bigskip

\noindent A.5. {\bf $D_\el$ case.}
We have:
\begin{align*}
\La_i(z) &= Y_i(zq^{-i+1}) Y_{i-1}(zq^{-i})^{-1}, \quad \quad
i=1,\ldots,\el-2, \\
\La_{\el-1}(z) &= Y_\el(zq^{-\el+2}) Y_{\el-1}(zq^{-\el+2})
Y_{\el-2}(zq^{-\el+1})^{-1},
\\
\La_\el(z) &= Y_{\el-1}(zq^{-\el+2}) Y_\el(zq^{-\el})^{-1}, \\
\La_{\ol{\el}}(z) &= Y_\el(zq^{-\el+2}) Y_{\el-1}(zq^{-\el})^{-1}, \\
\La_{\ol{\el-1}}(z) &= Y_{\el-2}(zq^{-\el+1}) Y_{\el-1}(zq^{-\el})^{-1}
Y_\el(zq^{-\el})^{-1}, \\
\La_{\ol{i}}(z) &= Y_{i-1}(zq^{-2\el+i+2}) Y_i(zq^{-2\el+i+1})^{-1}.
\end{align*}

Let
$$\s_i(z) = \sum_{\{ j_1,\ldots,j_i \} \in S} \La_{j_1}(z) \La_{j_2}(zq^2)
\ldots \La_{j_i}(zq^{2i-2}), \quad \quad i=1,\ldots,\el-2,$$ where $S$ is the
set of $\{ j_1,\ldots,j_i \}$, such that $j_\al<j_{\al+1}$ or
$j_\al=j_{\al+1}+1=\el+1, \al=1,\ldots,i-1$.

The formulas for $\s_{\el-1}(z)$ and $\s_\el(z)$, which correspond to the
spinor representations of $D_\el^{(1)}$ are given in \cite{FR:crit}.

The following formulas are due to A.~Kogan \cite{Ko}.

\begin{align*}
\{ \La_i(z),\La_i(w) \} &= {\mathcal M}_{11}\left( \frac{w}{z} \right)
\La_i(z)
\La_i(w), \\
\{ \La_i(z),\La_j(w) \} &= {\mathcal M}_{11}\left( \frac{w}{z} \right)
\La_i(z)
\La_j(w) +
\left( \delta \left( \frac{w}{zq^2} \right) - \delta \left( \frac{w}{z}
\right) \right. \\ &+ \left. \delta \left( \frac{w}{zq^{2\el-2i}} \right)
\delta_{i,\ol{j}} - \delta \left( \frac{w}{zq^{2\el-2i-2}} \right)
\delta_{i,\ol{j}} \right) \La_i(z) \La_j(w),
\end{align*}
if $i<j$. 

\begin{align*}
\{ T_1(z),T_1(w) \} &= {\mathcal M}_{11} \left( \frac{w}{z} \right) T_1(z)
T_1(w) \\ &+ \delta \left( \frac{w}{zq^2} \right) T_2(z) - \delta \left(
\frac{wq^2}{z} \right) T_2(w) \\ &+ \delta \left( \frac{w}{zq^{2\el-2}}
\right) - \delta \left( \frac{w^{2\el-2}}{z} \right).
\end{align*}

\section*{Appendix B. The difference version of the Drinfeld-Sokolov
reduction in the case $\g=C_2$}

\subsection*{B.1. Recollections}

Let us briefly recall the difference Drinfeld-Sokolov reduction scheme
\cite{FRS,SS}. Let $G$ be a simple algebraic group, and $\GG=G((z))$
be the corresponding formal loop group. Consider the following action
of $\GG$ on itself by difference gauge transformations:
\begin{equation}    \tag{B.1}
g(z) \cdot x(z) = g(zp) x(z) g(z)^{-1},
\end{equation}
where $p$ is a non-zero complex number. Let $s_i, i=1,\ldots,\el$, be
the simple reflections from the Weyl group $W=N(H)/H$, where $H$ is
the Cartan subgroup of $G$, and $N(H)$ is its normalizer. Choose once
and for all their lifts $n_i, i=1,\ldots,\el$, to $N(H) \subset
G$. Denote by $U$ the upper unipotent subgroup of $G$, and by $U_i$
its one-dimensional unipotent subgroup corresponding to the $i$th
simple root. Now set
$$
M=U n_1 \ldots n_\el U, \quad \quad
C = U_1 n_1 U_1 n_2 \ldots U_\el n_\el.
$$
Let $\UU=U((z)), \MM=M((z)), \CC=C((z))$.

\begin{thm}[\cite{SS}]
{\em The restriction of the action (B.1) to $\UU \subset \GG$
preserves $\MM$. The resulting action of $\UU$ on $\MM$ is free, and
$\CC$ is the cross-section of this action.}
\end{thm}

Furthermore, in \cite{FRS} (for $\g=A_1$) and \cite{SS} (in general) a
Poisson structure was defined on $\GG$, which descends to a
well-defined Poisson structure on $\MM/\UU \simeq \CC$. The
corresponding Poisson algebra of functions on $\CC$ is the Poisson
algebra $\W^p(\g)$ considered in Sect.~4.3.

Let ${\mathcal H}$ be the formal loop group of $H$. We define
following \cite{SS} the difference Miura transformation $m: {\mathcal
H} \arr \CC$. Let $\ol{U}$ be the lower unipotent subgroup of
$G$. There is a unique element $f \in M \cap \ol{U}$. Denote by $s$
the Coxeter element of $W$: $s=s_1 \ldots s_\el$ (note that $s$ here
corresponds to $s^{-1}$ of \cite{SS}). Define an embedding $i:
{\mathcal H} \arr \MM$ that sends $x \in {\mathcal H}$ to $x \cdot f
\cdot s^{-1}(x^{-1})$. It is easy to see that $i(x) \in \MM$.

The difference Miura transformation $m: {\mathcal H} \arr \CC$ is by
definition the composition of $i$ and the projection $\MM \arr
\CC$. In \cite{SS} a Poisson structure was defined on ${\mathcal H}$,
that makes the map $m$ Poisson.

The above construction has been worked out explicitly in the case
$\g=A_\el$ in \cite{FRS,SS}: it was shown that the Poisson algebra
$\W^t(A_\el)$ is isomorphic to the Poisson algebra $\W_{t,1}(A_\el)$,
and the Miura transformations for the two Poisson algebras
coincide. This proves \conjref{tds} in the case $\g=A_\el$. In the
next subsection we prove \conjref{tds} in the case $\g=C_2$.

\subsection*{B.2. The case of $C_2$}

We realize $C_2 = \on{Sp}_4$ as the group of $4 \times 4$ matrices
$A$ with $\det A = 1$, which satisfy
$$
A^t J A = J,
$$
where
$$
J = \begin{pmatrix} 0 & 0 & 0 & -1 \\
                    0 & 0 & -1 & 0 \\
                    0 & 1 & 0 & 0 \\
                    1 & 0 & 0 & 0 \\
\end{pmatrix}.
$$

We choose
$$
n_1 = \begin{pmatrix} 1 & 0 & 0 & 0 \\
                    0 & 0 & 1 & 0 \\
                    0 & -1 & 0 & 0 \\
                    0 & 0 & 0 & 1 \\
\end{pmatrix}, \quad \quad
n_2 = \begin{pmatrix} 0 & 1 & 0 & 0 \\
                    -1 & 0 & 0 & 0 \\
                    0 & 0 & 0 & -1 \\
                    0 & 0 & 1 & 0 \\
\end{pmatrix}.
$$

We have
$$
U = \left\{ \begin{pmatrix} 1 & a & b+\frac{ad}{2} & c \\
                    0 & 1 & d & b-\frac{ad}{2} \\
                    0 & 0 & 1 & -a \\
                    0 & 0 & 0 & 1 \\
\end{pmatrix}, a,b,c,d \in \C \right\},
$$
$$
U_1 = \left\{ \begin{pmatrix} 1 & 0 & 0 & 0 \\
                    0 & 1 & \al & 0 \\
                    0 & 0 & 1 & 0 \\
                    0 & 0 & 0 & 1 \\
\end{pmatrix} \right\}, \quad \quad
U_2 = \left\{ \begin{pmatrix} 1 & \beta & 0 & 0 \\
                    0 & 1 & 0 & 0 \\
                    0 & 0 & 1 & -\beta \\
                    0 & 0 & 0 & 1 \\
\end{pmatrix} \right\},
$$
and so
\begin{equation}    \tag{B.2}
C = \left\{ \begin{pmatrix} S_1 & 1 & 0 & 0 \\
                    -S_2 & 0 & S_1 & -1 \\
                    1 & 0 & 0 & 0 \\
                    0 & 0 & 1 & 0 \\
\end{pmatrix}, S_1, S_2 \in \C \right\}.
\end{equation}
Next, we choose:
$$
f = \begin{pmatrix} 1 & 0 & 0 & 0 \\
                    -1 & 1 & 0 & 0 \\
                    1 & -1 & 1 & 0 \\
                    1 & 0 & 1 & 1 \\
\end{pmatrix},
$$
and write a typical element of $H$ as
$$
x = \begin{pmatrix} x_2^{1/2} & 0 & 0 & 0 \\
                    0 & x_1 x_2^{-1/2}& 0 & 0 \\
                    0 & 0 & x_1^{-1} x_2^{1/2} & 0 \\
                    0 & 0 & 0 & x_2^{-1/2} \\
\end{pmatrix}.
$$

Now, to find explicit formulas for the difference Miura
transformation, we have to represent the matrix (B.2) as
$$
u(zt) \cdot x(z) \cdot f \cdot s^{-1}(x(z)^{-1}) \cdot u(z)^{-1},
$$
where $u(z) \in \UU$. Straightforward but tedious calculation gives:
\begin{align}    \tag{B.3}
S_1(z) &= x_1(z) + x_1(zp)^{-1} x_1(z)^2 x_2(z)^{-1} + x_1(zp)^{-1}
x_2(zp) \\ \notag &+ x_1(zp^2)^{-1}, \\    \tag{B.4}
S_2(z) &= x_1(z)^2 x_2(z)^{-1} + x_2(zp) + 2 x_1(zp) x_1(zp^2)^{-1}
\\ \notag &+ x_1(zp)^2 x_1(zp^2)^{-2} x_2(zp)^{-1} + x_1(zp^2)^{-2}
x_2(zp^2).
\end{align}

Now set $p=t^2$. If we let $Y_1(z) = x_1(z)$, and $Y_2(z) =
x_1(zt^{-1})^2 x_2(zt^{-1})^{-1}$, then formula (B.3) becomes:
\begin{equation}    \tag{B.5}
S_1(z) = Y_1(z) + Y_2(zt) Y_1(zt^2)^{-1} + Y_2(zt^3)^{-1} Y_1(zt^2) +
Y_1(zt^4)^{-1}.
\end{equation}
It coincides with the formula for the generator $T_1(z)$ of
$\W_{1,t}(C_2)$.

On the other hand, if we let $\wt{Y}_1(z) = x_1(z)^2 x_2(z)^{-1}$, and
$\wt{Y}_2(z) = x_1(zt)$, then formula (B.4) becomes:
\begin{align}    \tag{B.6}
S_2(z) &= \wt{Y}_1(z) + \wt{Y}_2(zt)^2 \wt{Y}_1(zt^2)^{-1} + 2
\wt{Y}_2(zt) \wt{Y}_2(zt^3)^{-1} \\ \notag &+ \wt{Y}_2(zt^3)^{-2}
\wt{Y}_1(zt^2) + \wt{Y}_1(zt^4)^{-2}.
\end{align}
It coincides with the formula for the generator $T_1(z)$ of
$\W_{1,t}(B_2)$, which is the same as the generator $T_2(z)$ of
$\W_{1,t}(C_2)$ up to a shift. Thus, the generators of $\W_{1,t}(C_2)$
and $\W^{t^2}(C_2)$ coincide.

Now we compute the Poisson brackets between $Y_i(z)$ using formula
(4.1) of \cite{SS}. First we find the eigenvectors of the
transformation $x \arr s^{-1}(x)$ on $H$. These are represented by
matrices
$$
e_1 = \begin{pmatrix} 1 & 0 & 0 & 0 \\
                       0 & i & 0 & 0 \\
                       0 & 0 & -i & 0 \\
                       0 & 0 & 0 & 1
\end{pmatrix}, \quad \quad
e_2 = \begin{pmatrix} 1 & 0 & 0 & 0 \\
                       0 & -i & 0 & 0 \\
                       0 & 0 & i & 0 \\
                       0 & 0 & 0 & 1
\end{pmatrix},
$$
with the eigenvalues $-i$ and $i$, respectively. In terms of these
elements, the Poisson tensor of \cite{SS} on ${\mathcal H}$ is
represented by the formal power series
\begin{equation}    \tag{B.7}
- \sum_{n\in\Z} \left( \frac{w}{z} \right)^n (1-t^{2n}) \frac{1}{2}
 \left( \frac{1-i}{1-it^{2n}} e_1 \otimes e_2 + \frac{1+i}{1+it^{2n}}
 e_2 \otimes e_1 \right).
\end{equation}
(recall that $p=t^2$).

Now let $B_i(z)$ be $(i,i)$th matrix entry of element of ${\mathcal
H}$. Using formula (B.7), we obtain the following Poisson brackets:
\begin{align*}
\{ B_1(z),B_1(w) \} &= \sum_{n \in \Z} \left( \frac{w}{z} \right)^n
(1-t^{2m}) \frac{1}{2} \left( \frac{1-i}{1-it^{2n}} +
\frac{1+i}{1+it^{2n}} \right) B_1(z) B_1(w) \\ &= \sum_{n \in \Z}
\left( \frac{w}{z} \right)^n \frac{t^{2n} - t^{-2n}}{t^{2n} + t^{-2n}}
B_1(z) B_1(w),
\end{align*}
and similarly:
\begin{align*}
\{ B_1(z),B_2(w) \} &= \sum_{n \in \Z} \left( \frac{w}{z} \right)^n
\frac{(t^{n} - t^{-n})^2}{t^{2n} + t^{-2n}} B_1(z) B_2(w), \\
\{ B_2(z),B_2(w) \} &= \sum_{n \in \Z} \left( \frac{w}{z} \right)^n
\frac{t^{2n} - t^{-2n}}{t^{2n} + t^{-2n}} B_2(z) B_2(w),
\end{align*}

Under the map $x \arr x \cdot s^{-1}(x^{-1})$, $x_1(z)$ goes to
$B_1(z)$, and $x_1(z) x_2(z)^{-1}$ goes to $B_2(z)$. Hence $Y_1(z) =
B_1(z), Y_2(z) = B_1(zt^{-1}) B_2(zt^{-1})$, and straightforward
calculation gives us the following Poisson brackets:
\begin{align*}
\{ Y_1(z),Y_1(w) \} &= \sum_{n \in \Z} \left( \frac{w}{z} \right)^n
\frac{t^{2n} - t^{-2n}}{t^{2n} + t^{-2n}} Y_1(z) Y_1(w), \\
\{ Y_1(z),Y_2(w) \} &= \sum_{n \in \Z} \left( \frac{w}{z} \right)^n
\frac{2(t^{n} - t^{-n})}{t^{2n} + t^{-2n}} Y_1(z) Y_2(w), \\
\{ Y_2(z),Y_2(w) \} &= \sum_{n \in \Z} \left( \frac{w}{z} \right)^n
\frac{2(t^{2n} - t^{-2n})}{t^{2n} + t^{-2n}} Y_2(z) Y_2(w).
\end{align*}
This agrees with the Poisson bracket
$$
\{ Y_i(z),Y_j(w) \} = \sum_{n \in \Z} \left( \frac{w}{z} \right)^n
(t^n-t^{-n}) M_{C_2}(1,t^n) \; Y_i(z) Y_j(w),
$$
in ${\mathcal H}_{1,t}(C_2)$. Hence $\W_{1,t}(C_2)$ and
$\W^{t^2}(C_2)$ are isomorphic as Poisson algebras, which is what we
wanted to show.

\section*{Appendix C. The matrices $M(q,t)$ for Lie algebras of
classical types}

\noindent {\bf $A_\el$ series.}

\begin{equation*}
M_{ij}(q,t) =
\frac{(q^it^{-i}-q^{-i}t^i)(q^{\el+1-j}t^{-\el-1+j}-q^{-\el-1+j}
t^{\el+1-j})}{q^{\el+1}t^{-\el-1}-q^{-\el-1}t^{\el+1}}, \quad
\quad i\leq j.
\end{equation*}

\bigskip

\noindent {\bf $B_\el$ series.}

\begin{align*}
M_{ij}(q,t) &=
\frac{(q^{2i} t^{-i} - q^{-2i} t^i)(q^{2\el-1-2j} t^{-\el+j} +
q^{-2\el+1+2j} t^{\el-i})(q+q^{-1})}{(q^2 t^{-1} - q^{-2} t)(q^{2\el-1}
t^{-\el}+ q^{-2\el+1} t^\el)}, & 1\leq i\leq j<\el, \\
M_{i\el}(q,t) &= \frac{(q^{2i} t^{-i} -q^{-2i} t^i)(q+q^{-1})}{(q^2
t^{-1} - q^{-2} t)(q^{2\el-1} t^{-\el}+ q^{-2\el+1} t^\el)}, &
1\leq i<\el, \\
M_{\el\el}(q,t) &=
\frac{q^{2\el} t^{-\el}-q^{-2\el} t^\el}{(q^2 t^{-1} - q^{-2}
t)(q^{2\el-1} t^{-\el}+ q^{-2\el+1} t^\el)}. &
\end{align*}

\bigskip

\noindent {\bf $C_\el$ series.}

\begin{equation*}
M_{ij}(q,t) =
\frac{(q^i t^{-i}-q^{-i}t^i)(q^{\el+1-j} t^{-\el+j} + q^{-\el-1+j}
t^{\el-j})}{(qt^{-1} - q^{-1} t)(q^{\el+1} t^{-\el} + q^{-\el-1}
t^\el)}, \quad \quad 1\leq i\leq j\leq \el.
\end{equation*}

\bigskip

\noindent {\bf $D_\el$ series.}

\begin{align*}
\displaystyle M_{ij}(q,t) &= \displaystyle \frac{(q^i t^{-i} -
q^{-i} t^i)(q^{\el - 1 - j} t^{-(\el - 1 - j)} + q^{-(\el - 1 - j)}
t^{\el - 1 - j})}{\den}, & 1\leq i\leq j\leq \el-1,\\
\displaystyle M_{i\el}(q,t) &= \displaystyle
\frac{q^i t^{-i} - q^{-i}t^i}{\den}, & 1\leq i\leq \el-2, \\ 
\displaystyle M_{\el-1,\el}(q,t) &= \displaystyle \frac{q^{\el-2}t^{-\el +
2} - q^{-\el +2}t^{\el-2}}{(\con)(\den)}, & \\
\displaystyle M_{\el-1,\el-1}(q,t) &= M_{\el\el}(q,t) = \displaystyle
\frac{q^{\el} t^{-\el} - q^{-\el} t^\el}{(\con)(\den)}. &
\end{align*}

\section*{Appendix D. Deformed Chiral Algebras}

\noindent {\bf Definition.} A deformed chiral algebra (DCA) with a
representation is a collection of the following data:
\begin{itemize}
\item A vector space $V$ called the space of fields.

\medskip
\item A vector space $W = \cup_{n\geq 0} W_n$ called the space of states,
which is union of finite-di\-men\-sional subspaces $W_n$. We consider a
topology on $W$ in which $\{ W_n \}_{n\geq 0}$ is the base of open
neighborhoods of $0$.

\medskip
\item A linear map $Y: V\to \on{End} W \widehat{\otimes} [[z, z^{-1}]]$
that is for each $A \in V$, a formal power series $Y(A,z) = \sum_{n\in\Z}
A_n z^{-n}$, where each $A_n$ is a linear operator $W \arr W$, such that
$A_n \cdot W_m \subset W_{m+N(n)}, \forall m \geq 0$, for some $N(n) \in
\Z$.

\medskip
\item A meromorphic function $S(x): \C^\times \arr \on{Aut}(V \otimes V)$,
satisfying the Yang-Baxter equation:
$$S_{12}(z) S_{13}(zw) S_{23}(w) =  S_{23}(w) S_{13}(zw) S_{12}(z),
$$ for
all $z,w \in \C^\times$.

\medskip
\item A lattice $L \subset \C^\times$, which contains the poles of
$S(x)$.

\medskip
\item A vector $\Omega \in V$, such that $Y(\Omega,z) = \on{Id}$.
\end{itemize}

The data of deformed chiral algebra should satisfy the following axioms:
\begin{itemize}
\item[(1)] For any $A_i \in V, i=1,\ldots,n$, the composition $Y(A_1,z_1)
\ldots Y(A_n,z_n)$ converges in the domain $|z_1| \gg \ldots \gg |z_n|$ and
can be continued to a meromorphic operator valued function
$$R(Y(A_1,z_1)
\ldots Y(A_n,z_n)): (\C^\times)^n \arr \on{Hom}(W,\overline{W}),
$$ where
$\overline{W}$ is the completion of $W$ with respect to its topology.

\medskip
\item[(2)] Denote $R(Y(A,z)Y(B,w))$ by $Y(A \otimes B;z,w)$. Then
$$Y(A \otimes B;z,w) = Y(S(w/z)(B \otimes A);w,z).$$

\medskip
\item[(3)] The poles of the meromorphic function $R(Y(A,z)Y(B,w))$ lie on
the lines $z = w\gamma$ where $\gamma\in L$. For each such line and $n\geq
0$, there exists $C_n \in V$, such that
$$
\Res_{z=w\gamma} R(Y(A,z)Y(B,w)) (z-w\gamma)^n \dzz = Y(C_n,w).
$$
\end{itemize}

\end{document}